\documentclass[12pt]{amsart}
\usepackage[francais,english]{babel}

\usepackage[]{amsfonts}
\usepackage{amssymb}
\usepackage{amsmath}

\usepackage[T1]{fontenc} %Si T1 ->fontes ec* qui tolerent les caracteres accentues sinon OT1

\def\couleur(#1 #2 #3)
	{
\special{ps::[end]}}

\def\bx#1{\setbox1=\hbox{\kern3pt{#1}\kern3pt}			% Make a box. Close it by "}"
 \dimen1=\ht1 \advance\dimen1 by 3pt \dimen2=\dp1 \advance\dimen2 by 3pt
 \setbox1=\hbox{\vrule height\dimen1 depth\dimen2\box1\vrule}%
 \setbox1=\vbox{\hrule\box1\hrule}%
 \advance\dimen1 by .4pt \ht1=\dimen1
 \advance\dimen2 by .4pt \dp1=\dimen2 \box1\relax}

\def\wbb#1{\kern#1em}
\def\vci{\vrule  width.02em height1.47ex depth-.0ex}		% le 1 en blackboard
\def\11{{\rm\wbb{.2}\vci\wbb{-.37}1}}

\newtheorem{Theorem}{Theorem}[section]
\newtheorem{Lemma}[Theorem]{Lemma}
\newtheorem{Definition}[Theorem]{Definition}
\newtheorem{Proposition}[Theorem]{Proposition}
\newtheorem{Corollary}[Theorem]{Corollary}
\newtheorem{Remark}[Theorem]{Remark}

\begin{document}

\title{A subordination principle. Applications.}

\author{Eric Amar}

\address{Universit\'e Bordeaux, IMB, UMR 5251, F-33400 Talence, France.}

\email{Eric.Amar@math.u-bordeaux1.fr}

\subjclass[2010]{32A50, 42B30}
\maketitle
 \ \par 
\ \par 

\renewcommand{\abstractname}{R\'esum\'e}

\begin{abstract}
Ce principe de subordination dit grossi\`erement : si une propri\'et\'e
 est vrai pour les espaces de Hardy pour certains domaines de
  ${\mathbb{C}}^{n}$  alors elle vrai pour les espaces de Bergman
 pour les domaines du m\^eme genre de  ${\mathbb{C}}^{n-1}.$ \ \par 
On donne des applications de ce principe aux mesures de Bergman-Carleson,
 aux suites d'interpolation pour les espaces de Bergman, au th\'eor\`eme
 de la couronne  $A^{p}$  et \`a la caract\'erisation des z\'eros
 de la classe de Bergman-Nevanlinna.\ \par 
Ces applications donnent des r\'esultats pr\'ecis pour les domaines
 born\'es strictement pseudo-convexes et les domaines born\'es
 convexes de type finis dans  ${\mathbb{C}}^{n}.$ \ \par 
\end{abstract}
\ \par 
\renewcommand{\abstractname}{Abstract}

\begin{abstract}
This subordination principle states roughly: if a property is
 true for Hardy spaces in some kind of domains in  ${\mathbb{C}}^{n}$
  then it is also true for the Bergman spaces of the same kind
 of domains in  ${\mathbb{C}}^{n-1}.$ \ \par 
We give applications of this principle to Bergman-Carleson measures,
 interpolating sequences for Bergman spaces,  $A^{p}$  Corona
 theorem and characterization of the zeros set of Bergman-Nevanlinna
 class.\ \par 
These applications give precise results for bounded strictly-pseudo
 convex domains and bounded convex domains of finite type in
  ${\mathbb{C}}^{n}.$ \ \par 
\end{abstract}
\ \par 
\ \par 

\tableofcontents
\ \par 

\section{Introduction.}
\quad  Let us start with some definitions. {\sl In all the sequel,
 } domain will mean bounded connected open set in  ${\mathbb{C}}^{n}$
  with smooth  ${\mathcal{C}}^{\infty }$  boundary defined by
 a real valued function  $r\in {\mathcal{C}}^{\infty }({\mathbb{C}}^{n}),$
 \ \par 
 i.e.             $\Omega =\lbrace z\in {\mathbb{C}}^{n}::r(z)<0\rbrace
 ,\ \forall z\in \partial \Omega ,\ {\rm{grad }}r(z)\neq 0,$ \ \par 
with the defining function  $r$  such that  $\forall z\in \Omega
 ,\ -r(z)\simeq d(z,\Omega ^{c})$  uniformly on  $\bar \Omega
 .$  (See the beginning of section~\ref{subPrinFin434} for the
 existence of such a function)\ \par 
\quad \quad  	Associate to it the "lifted" domain  $\tilde \Omega $  in 
 $(z,w)\in {\mathbb{C}}^{n+k}$  with defining function\ \par 
\quad \quad \quad  $\tilde r(z,w):=r(z)+\left\vert{w}\right\vert ^{2}.$ \ \par 
Usually our defining functions will be pluri-sub-harmonic, p.s.h.
 or even strictly pluri-sub-harmonic, s.p.s.h., in a neighborhood
 of  $\displaystyle \bar \Omega .$ \ \par 
This operation keeps the nature of the domain :\ \par 
\quad  $\bullet $  if  $\Omega $  is pseudo-convex defined by a  $r$
  p.s.h.,  $\tilde \Omega $  is still pseudo-convex defined by
  $\displaystyle \tilde r$  p.s.h.;\ \par 
\quad  $\bullet $  if  $\Omega $  is strictly pseudo-convex defined
 by a  $r$  s.p.s.h., so is  $\tilde \Omega \ ;$ \ \par 
\quad  $\bullet $  if  $\Omega $  is convex defined by a function 
 $r$  convex  , so is  $\tilde \Omega \ ;$ \ \par 
\quad  $\bullet $  if  $\Omega $  is convex of finite type  $m,$  defined
 by a function  $r$  convex, so is  $\tilde \Omega .$ \ \par 
Moreover we still have	 	 $\forall (z,w)\in \tilde \Omega ,\
 -(r(z)+\left\vert{w}\right\vert ^{2})\simeq d((z,w),\ \tilde
 \Omega ^{c}).$ \ \par 
\quad  Let  $\,dm(z)$  be the Lebesgue measure in  ${\mathbb{C}}^{n}$
  and  $\,d\sigma (z)$  be the Lebesgue measure on  $\partial \Omega .$ \ \par 
\quad  For  $z\in \Omega ,$  let  $\delta (z):=d(z,\Omega ^{c})\simeq
 -r(z)$  be the distance from  $z$  to the boundary of  $\Omega .$ \ \par 
\quad \quad  	For  $k\in {\mathbb{N}},$  let  $v_{k}$  be the volume of the
 unit ball in  ${\mathbb{C}}^{k}$  and set\ \par 
\quad \quad \quad \quad \quad 	 $\forall z\in \Omega ,\ dm_{0}(z):=dm(z),$ \ \par 
\quad \quad \quad \quad \quad   $\forall k\geq 1,\ \forall z\in \Omega ,\ dm_{k}(z):=(k+1)v_{k+1}(-r(z))^{k}dm(z)$
 \ \par 
a weighted Lebesgue measure in  $\Omega $  suitable for our
 needs. Clearly we have that  $dm_{k}(z)\simeq \delta (z)^{k}dm(z).$ \ \par 
\quad  Let  ${\mathcal{U}}$  be a neighbourhood of  $\partial \Omega
 $  in  $\Omega $  such that the normal projection  $\pi $  onto
  $\partial \Omega $  is a smooth well defined application.\ \par 
\quad  Define the Bergman, Hardy and Nevanlinna spaces as usual :\ \par 
\begin{Definition}
 Let  $f$  be a holomorphic function in  $\Omega \ ;$  we say
 that  $f\in A_{k}^{p}(\Omega )$  if\par 
\quad \quad \quad  $\displaystyle \ {\left\Vert{f}\right\Vert}_{k,p}^{p}:=\int_{\Omega
 }{\left\vert{f(z)}\right\vert ^{p}\,dm_{k}(z)}<\infty .$ \par 
\quad  We say that  $f\in {\mathcal{N}}_{k}(\Omega )$  if\par 
\quad \quad \quad  $\displaystyle \ {\left\Vert{f}\right\Vert}_{{\mathcal{N}}_{k}}=\int_{\Omega
 }{\log ^{+}\left\vert{f(z)}\right\vert \,dm_{k}(z)}<\infty ,$ \par 
\quad  We say that  $f\in H^{p}(\Omega )$  if\par 
\quad \quad \quad  $\displaystyle \ {\left\Vert{f}\right\Vert}_{p}^{p}:=\sup _{\epsilon
 >0}\int_{\lbrace r(z)=-\epsilon \rbrace }{\left\vert{f(\pi (z))}\right\vert
 ^{p}\,d\sigma (z)}<\infty .$ \par 
\quad  Finally we say that  $f\in {\mathcal{N}}(\Omega )$  if\par 
\quad \quad \quad  $\displaystyle \ {\left\Vert{f}\right\Vert}_{{\mathcal{N}}}=\sup
 _{\epsilon >0}\int_{\lbrace r(z)=-\epsilon \rbrace }{\log ^{+}\left\vert{f(\pi
 (z))}\right\vert \,d\sigma (z)}<\infty .$ \par 
\end{Definition}
This is meaningful because, for  $\epsilon $  small enough,
 the set  $\lbrace r(z)=-\epsilon \rbrace $  is a smooth manifold
 in  $\Omega $  contained in  ${\mathcal{U}}.$ \ \par 
\quad  Now we can state our subordination lemma:\ \par 
\ \par 
\ \par 
\begin{Theorem}
 ~\label{subPrinCarl24}(Subordination lemma) Let  $\Omega $ 
 be a domain in  ${\mathbb{C}}^{n},\ \tilde \Omega $  its lift
 in  ${\mathbb{C}}^{n+k}$  and  $F(z,w)\in H^{p}(\tilde \Omega
 ),$  we have  $f(z):=F(z,0)\in A_{k-1}^{p}(\Omega )$  and  
 $\ {\left\Vert{f}\right\Vert}_{A_{k-1}^{p}(\Omega )}\lesssim
 {\left\Vert{F}\right\Vert}_{H^{p}(\tilde \Omega )};$ \par 
if  $F(z,w)\in {\mathcal{N}}(\tilde \Omega ),$  then  $f(z):=F(z,0)\in
 {\mathcal{N}}_{k-1}(\Omega )$  and  $\ {\left\Vert{f}\right\Vert}_{{\mathcal{N}}_{k-1}(\Omega
 )}\lesssim {\left\Vert{F}\right\Vert}_{{\mathcal{N}}(\tilde \Omega )}.$ \par 
 A function  $f,$  holomorphic in  $\Omega ,$  is in the Bergman
 space  $A_{k-1}^{p}(\Omega )$  (resp. in the Nevanlinna Bergman
 space  ${\mathcal{N}}_{k-1}(\Omega )$ ) if and only if the function
  $F(z,w):=f(z)$  is in the Hardy space  $H^{p}(\tilde \Omega
 )$  (resp. in the Nevanlinna class  ${\mathcal{N}}(\tilde \Omega
 )$ ) and we have  $\ {\left\Vert{f}\right\Vert}_{A_{k-1}^{p}}\simeq
 {\left\Vert{F}\right\Vert}_{H^{p}(\tilde \Omega )}$  (resp.
  $\ {\left\Vert{f}\right\Vert}_{{\mathcal{N}}_{k-1}(\Omega )}\simeq
 {\left\Vert{F}\right\Vert}_{{\mathcal{N}}(\tilde \Omega )}$ ).\par 
\end{Theorem}
\quad  In the section~\ref{subPrinFin434} we prove the subordination
 lemma as a consequence of a disintegration of Lebesgue measure.\ \par 
\quad  In the section~\ref{subPrinFin435} we introduce the notion of
 a "good" family of polydiscs  ${\mathcal{P}},$  directly inspired
 by the work of Catlin~\cite{Catlin84} and introduced in~\cite{AmAspc09}
 together with a homogeneous hypothesis,  $(Hg).$ \ \par 
 This notion allows us to define geometric Carleson measure,
 denoted as  $\Lambda (\Omega ),$  for Hardy spaces and denoted
 as  $\Lambda _{k}(\Omega ),$   for Bergman spaces and to put
 it in relation with the Carleson embedding theorem still for
 these two classes of spaces.\ \par 
\quad  In subsection~\ref{subPrinFin436} we apply the subordination
 lemma to get a Bergman-Carleson embedding theorem from a Hardy-Carleson
 embedding one.\ \par 
\quad  The bounded strictly pseudo-convex domains have Hardy-Carleson
 embedding property by a result of H\"ormander~\cite{HormPSH67},
 hence they have the Bergman-Carleson embedding property by this result.\ \par 
A direct application of it is the following\ \par 
\begin{Corollary}
A positive Borel measure  $\mu $  in a strictly pseudo-convex
 domain  $\Omega $  in  ${\mathbb{C}}^{n}$  verifies\par 
\quad \quad \quad \quad \quad 	 $\displaystyle \forall p\geq 1,\ \forall f\in A_{k-1}^{p}(\Omega
 ),\ \int_{\Omega }{\left\vert{f}\right\vert ^{p}\,d\mu }\lesssim
 {\left\Vert{f}\right\Vert}_{A_{k-1}^{p}(\Omega )}$ \par 
iff :\par 
\quad \quad \quad \quad 	 $\displaystyle \forall a\in \Omega ,\ \mu (P_{a}(2))\lesssim
 \delta (a)^{n+k},$ \par 
where  $P_{a}(2)$  is the polydisc of the good family  ${\mathcal{P}}$
  centered at  $a$  and of "radius"  $2.$ \par 
\end{Corollary}
\quad  This characterization was already proved by Cima and Mercer~\cite{CimaMercer95}
 even for the spaces  $A_{\alpha }^{p}(\Omega )$  with  $\alpha
 \geq 0.$  So, in the case where  $\alpha $  is an integer we
 recover their characterization.\ \par 
\quad \quad  	M. Abate and A. Saraco~\cite{AbateSaraco11} studied Carleson
 measures in strongly pseudo-convex domains but with a different
 point of view : instead of using the family of polydiscs to
 characterize them they use invariant balls.\ \par 
\quad \quad  	We have also a characterization for convex domains of finite
 type, as shown in subsection~\ref{subPrinFin434}.\ \par 
\begin{Theorem}
 Let  $\Omega $  be a  convex domain of finite type in  ${\mathbb{C}}^{n}$
  ; the measure  $\mu $  verifies\par 
(*)      	      	 $\displaystyle \exists p>1,\ \exists C_{p}>0,\
 \forall f\in A_{k-1}^{p}(\Omega ),\ \int_{\Omega }{\left\vert{f}\right\vert
 ^{p}\,d\mu }\leq C_{p}^{p}{\left\Vert{f}\right\Vert}_{A_{k-1}^{p}(\Omega
 )}^{p}$ \par 
iff :\par 
(**)      	      	 $\displaystyle \exists C>0::\forall a\in
 \Omega ,\ \mu (\Omega \cap P_{a}(2))\leq Cm_{k-1}(\Omega \cap P_{a}(2)).$ \par 
Hence if   $\mu $  verifies (*) for a  $p>1,$  it verifies (*)
 for all  $q>1.$ \par 
\end{Theorem}
\ \par 
\quad  Now let  $\Omega $  be a domain in  ${\mathbb{C}}^{n}.$  We
 say that the  $H^{p}$ -Corona theorem is true for  $\Omega $
  if we have :\ \par 
\quad \quad \quad  $\forall g_{1},...,\ g_{m}\in H^{\infty }(\Omega )::\forall
 z\in \Omega ,\ \sum_{j=1}^{m}{\left\vert{g_{j}(z)}\right\vert
 }\geq \delta >0$ \ \par 
then\ \par 
\quad \quad \quad  $\forall f\in H^{p}(\Omega ),\ \exists (f_{1},...,\ f_{m})\in
 (H^{p}(\Omega ))^{m}::f=\sum_{j=1}^{m}{f_{j}g_{j}}.$ \ \par 
In the same vein, we say that the  $\displaystyle A_{k-1}^{p}(\Omega
 )$ -Corona theorem is true for  $\Omega $  if we have :\ \par 
\begin{equation} \forall g_{1},...,\ g_{m}\in H^{\infty }(\Omega
 )::\forall z\in \Omega ,\ \sum_{j=1}^{m}{\left\vert{g_{j}(z)}\right\vert
 }\geq \delta >0\end{equation}\ \par 
then\ \par 
\quad \quad \quad  $\forall f\in A_{k-1}^{p}(\Omega ),\ \exists (f_{1},...,\ f_{m})\in
 (A_{k-1}^{p}(\Omega ))^{m}::f=\sum_{j=1}^{m}{f_{j}g_{j}}.$ \ \par 
\quad  In the subsection~\ref{subPrinFin437}, we apply again the subordination
 principle, because the  $H^{p}$  Corona theorem is true in these
 cases, to get:\ \par 
\begin{Corollary}
 We have the  $\displaystyle A_{k}^{p}(\Omega )$ -Corona theorem
  in the following cases :\par 
\quad  $\bullet $  with  $p=2$  if  $\Omega $  is a bounded weakly
 pseudo-convex domain in  ${\mathbb{C}}^{n};$ \par 
\quad  $\bullet $  with  $1<p<\infty $  if  $\Omega $  is a bounded
  strictly pseudo-convex domain in  ${\mathbb{C}}^{n}.$ \par 
\end{Corollary}
\quad  In section~\ref{subPrinFin438} we define and study the interpolating
 sequences in a domain  $\Omega .$  We also define the notion
 of dual bounded sequences in  $H^{p}(\Omega )$  and in  $A_{k}^{p}(\Omega
 ),$  and applying the subordination principle to the result
 we proved for  $H^{p}(\Omega )$  interpolating sequences~\cite{AmAspc09},
 we get the following theorem.\ \par 
\begin{Theorem}
 If  $\Omega $  is a strictly pseudo-convex domain, or a convex
 domain of finite type in  ${\mathbb{C}}^{n}$  and if  $S\subset
 \Omega $  is a dual bounded sequence of points in  $A_{k}^{p}(\Omega
 )$  then, for any  $q<p,\ S$  is  $A_{k}^{p}(\Omega )$  interpolating
 with the linear extension property, provided that  $p=\infty
 $  or  $p\leq 2.$ \par 
\end{Theorem}
In the unit ball of  ${\mathbb{C}}^{n},$  we have a better result :\ \par 
\begin{Theorem}
 If  ${\mathbb{B}}$  is the unit ball in  ${\mathbb{C}}^{n}$
  and if  $S\subset {\mathbb{B}}$  is a dual bounded sequence
 of points in  $A_{k}^{p}({\mathbb{B}})$  then, for any  $q<p,\
 S$  is  $A_{k}^{p}(\Omega )$  interpolating with the linear
 extension property.\par 
\end{Theorem}
\ \par 
\quad  Finally in the section~\ref{subPrinFin439} we study zeros set
 for Nevanlinna Bergman functions.\ \par 
\quad  Let  $\Omega $  be a domain in  ${\mathbb{C}}^{n}$  and  $u$
  a holomorphic function in  $\Omega .$  Set  $X:=\lbrace z\in
 \Omega ::u(z)=0\rbrace $  the zero set of  $u$  and  $\Theta
 _{X}:=\partial \bar \partial \log  \left\vert{u}\right\vert
 $  its associated  $(1,1)$  current of integration.\ \par 
\begin{Definition}
 A zero set  $X$  of a holomorphic function  $u$  in the domain
  $\Omega $  is in the Blaschke class,  $X\in {\mathcal{B}}(\Omega
 ),$  if there is a constant  $C>0$  such that\par 
 \[\displaystyle \forall \beta \in \Lambda _{n-1,\ n-1}^{\infty
 }(\bar \Omega ),\ \left\vert{\int_{\Omega }{(-r(z))\Theta _{X}\wedge
 \beta }}\right\vert \leq C{\left\Vert{\beta }\right\Vert}_{\infty },\] \par 
where  $\displaystyle \Lambda _{n-1,\ n-1}^{\infty }(\bar \Omega
 )$  is the space of  $(n-1,n-1)$  continuous form in  $\bar
 \Omega ,$  equipped with the sup norm of the coefficients.\par 
\end{Definition}
\quad  If  $u\in {\mathcal{N}}(\Omega )$  then it is well known~\cite{zeroSkoda}
 that  $X$  is in the Blaschke class of  $\Omega .$ \ \par 
\quad  We do the analogue for the Bergman spaces :\ \par 
\begin{Definition}
 A zero set  $X$  of a holomorphic function  $u$  in the domain
  $\Omega $  is in the Bergman-Blaschke class,  $X\in {\mathcal{B}}_{k}(\Omega
 ),$  if there is a constant  $C>0$  such that\par 
 \[\displaystyle \forall \beta \in \Lambda _{n-1,\ n-1}^{\infty
 }(\bar \Omega ),\ \left\vert{\int_{\Omega }{(-r(z))^{k+1}\Theta
 _{X}\wedge \beta }}\right\vert \leq C{\left\Vert{\beta }\right\Vert}_{\infty
 },\] \par 
where  $\displaystyle \Lambda _{n-1,\ n-1}^{\infty }(\bar \Omega
 )$  is the space of  $(n-1,n-1)$  continuous form in  $\bar
 \Omega ,$  equipped with the sup norm of the coefficients.\par 
\end{Definition}
\quad  If  $u\in {\mathcal{N}}_{k}(\Omega )$  then   $X$  is in the
 Bergman-Blaschke class of  $\Omega $  as can be seen again by
 use of the subordination lemma.\ \par 
\quad  Hence exactly as for the Corona theorem we can set the definitions :\ \par 
we say that the {\sl Blaschke characterization is true for}
  $\Omega $  if we have :\ \par 
\quad \quad \quad  $X\in {\mathcal{B}}(\Omega )\Rightarrow \exists u\in {\mathcal{N}}(\Omega
 )$  such that  $X=\lbrace z\in \Omega ::u(z)=0\rbrace .$ \ \par 
And the same for the Bergman spaces :\ \par 
we say that the {\sl Bergman-Blaschke characterization is true
 for}  $\Omega $  if we have :\ \par 
\quad \quad \quad  $X\in {\mathcal{B}}_{k}(\Omega )\Rightarrow \exists u\in {\mathcal{N}}_{k}(\Omega
 )$  such that  $X=\lbrace z\in \Omega ::u(z)=0\rbrace .$ \ \par 
\quad  We get, by use of the subordination lemma applied to the corresponding
 Nevanlinna Hardy results,\ \par 
\begin{Corollary}
 The Bergman-Blaschke characterization is true  in the following cases :\par 
\quad  $\bullet $  if  $\Omega $  is a  strictly pseudo-convex domain
 in  ${\mathbb{C}}^{n}\ ;$ \par 
\quad  $\bullet $  if  $\Omega $  is a convex domain of finite type
 in  ${\mathbb{C}}^{n}.$ \par 
\end{Corollary}
\ \par 
\quad  We stated and proved the subordination lemma for the ball in
  ${\mathbb{C}}^{n}$  in 1978~\cite{amBerg78}, and, since then,
 we gave seminars and conferences about it in the general situation.\ \par 
\quad \quad  	As we seen some applications to strictly pseudo-convex domains
 done here are already known. The applications to the convex
 domains of finite type are new.\ \par 
\quad  I am grateful to Marco Abate for an interesting discussion on
 Bergman-Carleson measures in january 2010.\ \par 

\section{The subordination lemma.~\label{subPrinFin434}}
\quad \quad  	Let  $\Omega :=\lbrace z\in {\mathbb{C}}^{n}::\rho (z)<0\rbrace
 ,\ \partial \rho (z)\neq 0$  on  $\displaystyle \partial \Omega
 $  with  $\rho \in {\mathcal{C}}^{2}(\bar \Omega ).$   Let\ \par 
\quad \quad \quad \quad \quad   $\tilde \Omega :=\lbrace (z,w)\in {\mathbb{C}}^{n}{\times}{\mathbb{C}}::\rho
 (z)+\left\vert{w}\right\vert ^{2}<0\rbrace $  \ \par 
be the lift of  $\displaystyle \Omega $  in  ${\mathbb{C}}^{n+1}.$
  We can always manage to have  $\displaystyle \ \left\vert{{\rm{grad
 }}\rho (z)}\right\vert =1$  for  $\displaystyle z\in \partial
 \Omega $  by the well known following lemma.\ \par 
\begin{Lemma}
~\label{subPrinFin1351}Let  $\Omega $  be a domain in  ${\mathbb{R}}^{n},$
  we can always choose a defining function  $s$  for  $\displaystyle
 \Omega $  such that  $\displaystyle \forall z\in \partial \Omega
 ,\ \left\vert{{\rm{grad }}s(z)}\right\vert =1.$ \par 
\end{Lemma}
\quad \quad  	Proof.\ \par 
Because  ${\rm{grad }}\ r(z)\neq 0$  on  $\partial \Omega ,$
  we take any smooth strictly positive extension  $\displaystyle
 h$  of  $\displaystyle \ \frac{1}{\left\vert{{\rm{grad }}r(z)}\right\vert
 }$  in  $\displaystyle \bar \Omega \ ;$  then set  $\displaystyle
 s(z)=h(z)r(z).$  We have that  $\displaystyle {\rm{grad }}s(z)=h{\rm{grad
 }}r(z)+r(z){\rm{grad }}h(z)=h{\rm{grad }}r$  on  $\displaystyle
 \partial \Omega ,$  hence  $\displaystyle \ \left\vert{{\rm{grad
 }}s}\right\vert =1$  on  $\displaystyle \partial \Omega .$ 
 Of course because  $\displaystyle h>0,$  we have that  $\Omega
 =\lbrace z\in {\mathbb{R}}^{n}::s<0\rbrace .$   $\blacksquare $ \ \par 
\begin{Lemma}
 ~\label{subPrinFin1352}Let  $\Omega $  be a domain in  ${\mathbb{R}}^{n},$
  defined by a function  $r\in {\mathcal{C}}^{\infty }$ , i.e.\par 
\quad \quad \quad  $\Omega :=\lbrace x\in {\mathbb{R}}^{n}::r(z)<0\rbrace ,\ \forall
 x\in \partial \Omega ,\ {\rm{grad }}\ r(x)=1.$ \par 
Then the Lebesgue measure  $\sigma $  on  $\partial \Omega $  is given by\par 
\quad \quad \quad  $\displaystyle \forall g\in {\mathcal{C}}(\partial \Omega ),\
 \int_{\partial \Omega }{g\,d\sigma }=\lim  _{\eta \rightarrow
 0}\frac{1}{\eta }\int_{\lbrace -\eta \leq r(x)<0\rbrace }{\tilde
 g(x)\,dm(x)},$ \par 
where  $\tilde g(x)$  is any continuous extension of   $g$ 
 near  $\partial \Omega .$ \par 
\end{Lemma}
\quad \quad  	Proof.\ \par 
Because  $\displaystyle \partial \Omega $  is a codimension
 one manifold,  $\displaystyle \forall x\in \partial \Omega ,\
 {\rm{grad }}\ r(x)=1$  then  $\lbrace x\in {\mathbb{R}}^{n}::-\eta
 \leq r(x)<0\rbrace $  is "half" a tube of thickness  $\eta $
  around  $\displaystyle \partial \Omega ,$  hence we can apply
 corollary 6.9.12 in~\cite{GeoDiff87} or the original work by
 H. Weyl~\cite{HWeyl39}.  $\blacksquare $ \ \par 
\begin{Lemma}
~\label{subPrinFin1349} Let  $\displaystyle \Omega $  be a domain
 in  ${\mathbb{C}}^{n}.$  There is a defining function  $\rho
 $  for  $\displaystyle \Omega $  and  $\displaystyle \delta
 >0$  such that   $\displaystyle \ \left\vert{{\rm{grad }}\rho
 (z)}\right\vert ^{2}-4\rho (z)\geq \min (4\delta ,1/4).$ \par 
\end{Lemma}
\quad \quad  	Proof.\ \par 
Take a defining function  $\rho $  such that  $\displaystyle
 \ \left\vert{{\rm{grad }}\rho }\right\vert =1$  on  $\displaystyle
 \partial \Omega .$  Then the set  $\displaystyle K:=\lbrace
 z\in \Omega ::\left\vert{{\rm{grad }}\rho (z)}\right\vert \leq
 1/2\rbrace $  is compact in  $\displaystyle \Omega $  because
  $\displaystyle \ \left\vert{{\rm{grad }}\rho (z)}\right\vert
 $  is continuous. On this set  $K$  we have  $\displaystyle
 -\rho (z)\geq \delta >0$  because  $\displaystyle \rho (z)<0$
  in  $\displaystyle \Omega $  by definition of  $\displaystyle
 \Omega ,$  hence  $\displaystyle \rho (z)$  attains its maximum
  $\displaystyle -\delta <0$  on the compact  $K.$ \ \par 
\quad \quad  	Then we have\ \par 
\quad \quad \quad \quad \quad 	 $\displaystyle \forall z\in \bar \Omega ,\ $  $\displaystyle
 \ \left\vert{{\rm{grad }}\rho (z)}\right\vert ^{2}-4\rho (z)\geq
 \min (4\delta ,1/4),$ \ \par 
because\ \par 
\quad \quad 	 $\bullet $  in  $K,\ -\rho (z)\geq \delta \Rightarrow -4\rho
 (z)+\left\vert{{\rm{grad }}\rho (z)}\right\vert ^{2}\geq -4\rho
 (z)\geq 4\delta \ ;$ \ \par 
\quad \quad 	 $\bullet $  outside  $K,\ \left\vert{{\rm{grad }}\rho (z)}\right\vert
 >1/2\Rightarrow \left\vert{{\rm{grad }}\rho (z)}\right\vert
 ^{2}>1/4\Rightarrow \left\vert{{\rm{grad }}\rho (z)}\right\vert
 ^{2}-4\rho (z)\geq 1/4.$ \ \par 
Which completes the proof.  $\blacksquare $ \ \par 
\ \par 
\quad \quad  	Now back to the lifted domain  $\displaystyle \tilde \Omega
 .$ 	 The boundary of  $\displaystyle \tilde \Omega $  is defined
 by  $\displaystyle \rho (z)+\left\vert{w}\right\vert ^{2}=0,$
  hence on  $\displaystyle \partial \tilde \Omega $  we have
  $\displaystyle \ \left\vert{w}\right\vert ^{2}=-\rho (z).$ \ \par 
\begin{Lemma}
~\label{subPrinFin1350}Let  $\displaystyle \Omega $  be a domain
 in  ${\mathbb{C}}^{n}.$  There is a defining function  $\rho
 $  for  $\displaystyle \Omega $  and  $\displaystyle \delta
 >0$  such that  $\displaystyle \ \left\vert{{\rm{grad }}(\rho
 (z)+\left\vert{w}\right\vert ^{2})}\right\vert \geq \min (2\delta ,1/2).$ \par 
\end{Lemma}
\quad  Proof.\ \par 
Let us compute\ \par 
\quad \quad \quad \quad \quad   $\displaystyle {\rm{grad }}(\rho (z)+\left\vert{w}\right\vert
 ^{2})=(\frac{\partial \rho }{\partial x_{1}},\frac{\partial
 \rho }{\partial y_{1}},..,\frac{\partial \rho }{\partial x_{n}},\frac{\partial
 \rho }{\partial y_{n}},2u,2v)\ ;$ \ \par 
where  $\displaystyle z_{j}=x_{j}+iy_{j}$  and  $\displaystyle
 w=u+iv.$  Hence\ \par 
\quad \quad \quad \quad \quad 	 $\displaystyle \ \left\vert{{\rm{grad }}(\rho (z)+\left\vert{w}\right\vert
 ^{2})}\right\vert ^{2}=\left\vert{{\rm{grad }}\rho (z)}\right\vert
 ^{2}+4\left\vert{w}\right\vert ^{2}.$ \ \par 
By lemma~\ref{subPrinFin1349} we get on  $\displaystyle \partial
 \tilde \Omega ,$  replacing  $\delta $  by  $\displaystyle \delta
 ^{2},$ \ \par 
\quad \quad \quad \quad \quad 	 $\displaystyle \ \left\vert{{\rm{grad }}(\rho (z)+\left\vert{w}\right\vert
 ^{2})}\right\vert ^{2}=\left\vert{{\rm{grad }}\rho (z)}\right\vert
 ^{2}-4\rho (z)\geq \min (4\delta ^{2},1/4).$ \ \par 
Taking square root we get the lemma.  $\blacksquare $ \ \par 
\ \par 
\quad \quad  	Then we have the main lemma of this section on the disintegration
 of the Lebesgue measure  $\displaystyle d\tilde \sigma $  on
  $\displaystyle \partial \tilde \Omega $  :\ \par 
\begin{Lemma}
 (Main lemma) For any continuous function  $g$  on  $\displaystyle
 \ {\overline{\tilde \Omega }}\ :$ \par 
\quad \quad \quad \quad \quad 	 $\displaystyle \ \int_{\partial \tilde \Omega }{g(z,w)d\tilde
 \sigma (z,w)}=\int_{\Omega }{{\sqrt{-\rho (z)+\frac{\left\vert{{\rm{grad
 }}\rho (z)}\right\vert ^{2}}{4}}}\lbrace \int_{\left\vert{w}\right\vert
 ^{2}=-\rho (z)}{g(z,w)d\left\vert{w}\right\vert }\rbrace dm(z)},$ \par 
where  $\displaystyle d\left\vert{w}\right\vert $  is the {\bf
 normalized} Lebesgue measure on the circle  $\displaystyle \
 \left\vert{w}\right\vert ^{2}=-\rho (z)$  and  $\displaystyle
 dm(z)$  is the Lebesgue measure on  ${\mathbb{C}}^{n}.$ \par 
\end{Lemma}
\quad \quad  	Proof.\ \par 
we want a defining function whose gradient has norm  $1$  on
 the boundary, hence we set\ \par 
\quad \quad \quad \quad \quad 	 $\displaystyle \forall (z,w)\in \partial \tilde \Omega ,\ h(z,w):=\frac{1}{\
 \left\vert{{\rm{grad }}(\rho (z)+\left\vert{w}\right\vert ^{2})}\right\vert
 },$ \ \par 
because we have by lemma~\ref{subPrinFin1350} that  $\displaystyle
 \ \left\vert{{\rm{grad }}(\rho (z)+\left\vert{w}\right\vert
 ^{2})}\right\vert \geq \min (2\delta ,1/2)$  on  $\displaystyle
 \partial \tilde \Omega \ ;$  by continuity we have  $\displaystyle
 \ \left\vert{{\rm{grad }}(\rho (z)+\left\vert{w}\right\vert
 ^{2})}\right\vert \geq \frac{1}{2}\min (2\delta ,1/2)$  in a
 neighborhood  $V$  of  $\displaystyle \partial \tilde \Omega \ ;$ \ \par 
 as in lemma~\ref{subPrinFin1351} we set\ \par 
\quad \quad \quad \quad \quad 	 $\displaystyle \tilde \rho (z,w):=\frac{\rho (z)+\left\vert{w}\right\vert
 ^{2}}{\ \left\vert{{\rm{grad }}(\rho (z)+\left\vert{w}\right\vert
 ^{2})}\right\vert },$  for  $\displaystyle (z,w)\in V$  and
 we extend it to  $\displaystyle \tilde \Omega .$ \ \par 
\quad \quad  	Then  $\displaystyle \ \left\vert{{\rm{grad }}\tilde \rho (z,w)}\right\vert
 =1$  on  $\displaystyle \partial \tilde \Omega .$ \ \par 
Fix  $\eta >0$  and set  $\displaystyle \tilde \Omega _{\eta
 }:=\lbrace (z,w)\in {\mathbb{C}}^{n}{\times}{\mathbb{C}}::\tilde
 \rho (z,w)<-\eta \rbrace \subset \tilde \Omega \ ;$  let  $\eta
 $  be small enough such that  $\displaystyle \tilde \Omega \backslash
 \tilde \Omega _{\eta }\subset V.$ \ \par 
the Lebesgue measure on the manifold  $\displaystyle \partial
 \tilde \Omega $   can be defined by lemma~\ref{subPrinFin1352}
 this way :\ \par 
\quad \quad \quad \quad \quad 	 $\displaystyle I:=\int_{\partial \tilde \Omega }{g(z,w)d\sigma
 (z,w)}=\lim _{\eta \rightarrow 0}\frac{1}{\eta }\int_{\tilde
 \Omega \backslash \tilde \Omega _{\eta }}{g(z,w)dm(z,w)}.$ \ \par 
Hence, by Fubini,\ \par 
\quad \quad \quad \quad \quad 	 $\displaystyle \ \int_{\tilde \Omega \backslash \tilde \Omega
 _{\eta }}{g(z,w)dm(z,w)}=\int_{\Omega }{\lbrace \int_{-\eta
 \leq \tilde \rho (z,w)<0}{g(z,w)dm(w)}\rbrace dm(z)}.$ \ \par 
Fix  $\displaystyle z\in \Omega $  and let us study\ \par 
\quad \quad \quad \quad \quad 	 $\displaystyle -\eta \leq \frac{\rho (z)+\left\vert{w}\right\vert
 ^{2}}{\ \left\vert{{\rm{grad }}(\rho (z)+\left\vert{w}\right\vert
 ^{2})}\right\vert }<0\Rightarrow \rho (z)+\left\vert{w}\right\vert
 ^{2}<0\Rightarrow \left\vert{w}\right\vert ^{2}<-\rho (z).$ \ \par 
Recall that\ \par 
\quad \quad \quad \quad \quad   $\displaystyle \ \left\vert{{\rm{grad }}(\rho (z)+\left\vert{w}\right\vert
 ^{2})}\right\vert ^{2}=\left\vert{{\rm{grad }}\rho (z)}\right\vert
 ^{2}+4\left\vert{w}\right\vert ^{2}\Rightarrow \ \left\vert{{\rm{grad
 }}(\rho (z)+\left\vert{w}\right\vert ^{2})}\right\vert ={\sqrt{\left\vert{{\rm{grad
 }}\rho (z)}\right\vert ^{2}+4\left\vert{w}\right\vert ^{2}}}.$ \ \par 
The other side of the inequality gives\ \par 
\quad \quad \quad \quad   $\displaystyle -\eta {\sqrt{\left\vert{{\rm{grad }}\rho (z)}\right\vert
 ^{2}+4\left\vert{w}\right\vert ^{2}}}\leq \rho (z)+\left\vert{w}\right\vert
 ^{2}<0,$ \ \par 
hence raising to the square\ \par 
\quad \quad \quad  $\displaystyle (\rho (z)+\left\vert{w}\right\vert ^{2})^{2}\leq
 \eta ^{2}(\left\vert{{\rm{grad }}\rho (z)}\right\vert ^{2}+4\left\vert{w}\right\vert
 ^{2}).$ 	\ \par 
Set\ \par 
\quad \quad \quad \quad \quad   $\displaystyle a:=-\rho (z)\geq 0,\ b:=\left\vert{{\rm{grad
 }}\rho (z)}\right\vert ^{2}>0,\ X:=\left\vert{w}\right\vert
 ^{2}\geq 0,$ \ \par 
then this inequality becomes\ \par 
\quad \quad \quad \quad \quad 	 $\displaystyle (X-a)^{2}\leq \eta ^{2}(b+4X)\Rightarrow X^{2}-2(a+2\eta
 ^{2})X+a^{2}-\eta ^{2}b\leq 0.$ \ \par 
\quad \quad  	This implies that  $X$  must be between the  $2$  roots :\ \par 
\quad \quad \quad \quad \quad 	 $\displaystyle \Delta ^{2}:=(a+2\eta ^{2})^{2}-(a^{2}-\eta
 ^{2}b)=\eta ^{2}(4a+b+4\eta ^{2})\ ;$ \ \par 
hence the roots are\ \par 
\quad \quad \quad \quad \quad 	 $\displaystyle X':=(a+2\eta ^{2})-\eta {\sqrt{4a+b+4\eta ^{2}}}\
 ;\ X'':=(a+2\eta ^{2})+\eta {\sqrt{4a+b+4\eta ^{2}}}.$ \ \par 
We already have that  $\displaystyle \ \left\vert{w}\right\vert
 ^{2}=X<a=-\rho (z),$  hence, setting  $\displaystyle c(\eta
 ):=(a+2\eta ^{2})-\eta {\sqrt{4a+b+4\eta ^{2}}}.$ \ \par 
\quad \quad \quad \quad   $\displaystyle -\eta \leq \tilde \rho (z,w)<0\iff c(\eta )\leq
 \left\vert{w}\right\vert ^{2}<a.$ 	 $\displaystyle $ \ \par 
\ \par 
Now,  $g$  being continuous on  $\displaystyle \ {\overline{\tilde
 \Omega }},$  we get, with  $\displaystyle w=re^{i\theta }$ 
 in polar coordinates,\ \par 
\quad \quad \quad \quad \quad 	 $g(z,re^{i\theta })=g(z,{\sqrt{-\rho (z)}}e^{i\theta })+i(\eta ),$ \ \par 
the  $i(\eta )$  being uniform with respect to  $\displaystyle
 z,w$  in  $\displaystyle V.$  So let\ \par 
\quad \quad \quad \quad \quad 	 $\displaystyle J:=\frac{1}{\eta }\int_{-\eta \leq \tilde \rho
 (z,w)<0}{g(z,w)dm(w)}$ \ \par 
we have\ \par 
\quad \quad 	 $\displaystyle J=\frac{1}{\eta }\int_{c(\eta )\leq \left\vert{w}\right\vert
 ^{2}<a}{g(z,w)dm(w)}\ ;$ \ \par 
computing with polar coordinates,\ \par 
\quad \quad \quad \quad \quad 	 $\displaystyle J=\frac{1}{\eta }\int_{{\sqrt{c(\eta )}}}^{{\sqrt{a}}}{\lbrace
 \int_{0}^{2\pi }{(g(z,{\sqrt{-\rho (z)}}e^{i\theta })+i(\eta
 ))\frac{d\theta }{2\pi }}\rbrace rdr},$ \ \par 
hence\ \par 
\quad \quad \quad \quad \quad 	 $\displaystyle J=\int_{0}^{2\pi }{(g(z,{\sqrt{-\rho (z)}}e^{i\theta
 })+i(\eta ))\frac{d\theta }{2\pi }}{\times}\frac{1}{\eta }\int_{{\sqrt{c(\eta
 )}}}^{{\sqrt{a}}}{rdr},$ \ \par 
but\ \par 
\quad \quad \quad \quad   $\displaystyle \ \frac{1}{\eta }\int_{{\sqrt{c(\eta )}}}^{{\sqrt{a}}}{rdr}=\frac{1}{2\eta
 }(a-c(\eta ))=\frac{1}{2\eta }(a\ -((a+2\eta ^{2})-\eta {\sqrt{4a+b+4\eta
 ^{2}}}))={\sqrt{a+\frac{b}{4}+\eta ^{2}}}-\eta .$ \ \par 
so we get\ \par 
\quad \quad \quad \quad \quad 	 $\displaystyle J=({\sqrt{a+\frac{b}{4}+\eta ^{2}}}-\eta )\int_{0}^{2\pi
 }{(g(z,{\sqrt{-\rho (z)}}e^{i\theta })+i(\eta ))\frac{d\theta
 }{2\pi }}.$ \ \par 
Hence, letting  $\displaystyle \eta \rightarrow 0,$  we get\ \par 
\quad \quad \quad \quad \quad 	 $\displaystyle J\rightarrow {\sqrt{a+\frac{b}{4}}}\int_{0}^{2\pi
 }{g(z,{\sqrt{-\rho (z)}}e^{i\theta })\frac{d\theta }{2\pi }}.$ \ \par 
Putting it in  $I$ \ \par 
\quad \quad \quad \quad \quad 	 $\displaystyle I=\int_{\Omega }{{\sqrt{a+\frac{b}{4}}}\lbrace
 \int_{0}^{2\pi }{g(z,{\sqrt{-\rho (z)}}e^{i\theta })\frac{d\theta
 }{2\pi }}\rbrace dm(z)},$ \ \par 
i.e.\ \par 
\quad \quad \quad \quad \quad 	 $\displaystyle I=\int_{\Omega }{{\sqrt{-\rho (z)+\frac{\left\vert{{\rm{grad
 }}\rho (z)}\right\vert ^{2}}{4}}}\lbrace \int_{\left\vert{w}\right\vert
 ^{2}=-\rho (z)}{g(z,w)d\left\vert{w}\right\vert }\rbrace dm(z)},$ \ \par 
with  $\displaystyle d\left\vert{w}\right\vert $  the normalized
 Lebesgue measure on the circle 	 $\displaystyle \lbrace \left\vert{w}\right\vert
 ^{2}=-\rho (z)\rbrace .$   $\blacksquare $ \ \par 
\begin{Corollary}
 ~\label{subPrinFin1453}Setting  $\displaystyle h(z):={\sqrt{-\rho
 (z)+\frac{\left\vert{{\rm{grad }}\rho (z)}\right\vert ^{2}}{4}}},$
  we have that  $\displaystyle \exists \alpha >0,\ \beta >0$  such that\par 
\quad \quad 	 $\bullet $   $\displaystyle \forall z\in \bar \Omega ,\ \alpha
 \leq h(z)\leq \beta \ ;$ \par 
\quad \quad 	 $\bullet $    $\displaystyle \forall g\in {\mathcal{C}}(\partial
 \tilde \Omega ),\ \int_{\partial \tilde \Omega }{g(z,w)d\sigma
 (z,w)}=\int_{\Omega }{h(z)\lbrace \int_{\left\vert{w}\right\vert
 ^{2}=-\rho (z)}{g(z,w)d\left\vert{w}\right\vert }\rbrace dm(z)}.$ \par 
\quad \quad 	 $\bullet $    $\displaystyle \forall f\in {\mathcal{C}}(\partial
 \tilde \Omega ),\ \int_{\partial \tilde \Omega }{f(z,w)\frac{1}{h(z)}d\sigma
 (z,w)}=\int_{\Omega }{\lbrace \int_{\left\vert{w}\right\vert
 ^{2}=-\rho (z)}{f(z,w)d\left\vert{w}\right\vert }\rbrace dm(z)}.$ \par 
\end{Corollary}
\quad \quad  	Proof.\ \par 
We have  $\alpha =\min (\delta ,1/16)$  by lemma~\ref{subPrinFin1349}
 and  $\displaystyle \beta ={\left\Vert{h}\right\Vert}_{\infty
 }<\infty $  because  $h$  is continuous on  $\displaystyle \bar
 \Omega $  and  $\displaystyle \bar \Omega $  is compact.\ \par 
The second point is the main lemma.\ \par 
So it remains to prove the last assertion and for it we set
  $\displaystyle g(z,w):=\frac{f(z,w)}{h(z)}\in {\mathcal{C}}(\partial
 \tilde \Omega )$  and we apply the main lemma.  $\blacksquare $ \ \par 
\begin{Remark}
 In the case of the unit ball  ${\mathbb{B}}$  in  ${\mathbb{C}}^{n}$
  we get, with  $\displaystyle \rho (z)=\left\vert{z}\right\vert
 ^{2}-1$ as defining function, that  $\displaystyle -\rho (z)+\frac{\left\vert{{\rm{grad
 }}\rho (z)}\right\vert ^{2}}{4}=1,$  hence we have a disintegration
 of the Lebesgue measure on  $\partial \tilde {\mathbb{B}}$ 
 without weight.\par 
\end{Remark}
\ \par 
\quad  Now we can prove our subordination lemma~\ref{subPrinCarl24}
 stated in the introduction.\ \par 
We copy from~\cite{amBerg78}, and adapt from the ball to this
 general case. We shall prove it with several steps.\ \par 
\begin{Proposition}
 ~\label{subPrinFin1454}Let  $\displaystyle \Omega $  be a domain
 in  ${\mathbb{C}}^{n}$  and  $\displaystyle \tilde \Omega $
  its lift in  ${\mathbb{C}}^{n+1}.$  There are constants  $\alpha
 >0,\ \beta >0$  depending only on  $\displaystyle \Omega $ 
 such that if  $\displaystyle F\in $   $\displaystyle H^{p}(\tilde
 \Omega )$  then  $\displaystyle F(z,0)\in A^{p}(\Omega )$  and
  $\displaystyle \ {\left\Vert{F(\cdot ,0)}\right\Vert}_{A^{p}(\Omega
 )}\leq \frac{1}{\alpha }{\left\Vert{F}\right\Vert}_{H^{p}(\tilde
 \Omega )}.$ \par 
\quad \quad  	Conversely if  $\displaystyle f\in A^{p}(\Omega ),$  for  $\displaystyle
 (z,w)\in \tilde \Omega $  set  $\displaystyle F(z,w):=f(z),$
  then we have  $\displaystyle \ {\left\Vert{F}\right\Vert}_{H^{p}(\tilde
 \Omega )}\leq \beta {\left\Vert{f}\right\Vert}_{A^{p}(\Omega )}.$ \par 
\end{Proposition}
\quad \quad  	Proof.\ \par 
 If  $F(z,w)\in H^{p}(\tilde \Omega )$  we have\ \par 
\quad \quad \quad  $\displaystyle \ {\left\Vert{F}\right\Vert}_{p}^{p}:=\sup  _{\epsilon
 >0}\int_{\lbrace \tilde r(z,w)=-\epsilon \rbrace }{\left\vert{F(z,w)}\right\vert
 ^{p}\,d\tilde \sigma (z,w)}<\infty .$ \ \par 
Fix  $\epsilon >0$  and set  $\tilde \Omega =\tilde \Omega _{\epsilon
 }:=\lbrace (z,w)\in {\mathbb{C}}^{n+1}::\tilde r(z,w)<\epsilon
 \rbrace $  to apply what precede.\ \par 
By corollary~\ref{subPrinFin1453} the Lebesgue measure on  $\partial
 \tilde \Omega $  is\ \par 
\quad \quad \quad \quad \quad 	 $\displaystyle \forall g\in {\mathcal{C}}({\overline{\tilde
 \Omega }}),\ \int_{\partial \tilde \Omega }{g(z,w)d\tilde \sigma
 (z,w)}=\int_{\Omega }{h(z)\lbrace \int_{\left\vert{w}\right\vert
 ^{2}=-\rho (z)}{g(z,w)d\left\vert{w}\right\vert }\rbrace dm(z)},$ \ \par 
with  $\displaystyle \forall z\in \bar \Omega ,\ 0<\alpha \leq
 h(z)\leq \beta <\infty .$ \ \par 
So\ \par 
\quad \quad \quad \quad \quad 	 $\displaystyle \ \int_{\Omega }{h(z)\lbrace \int_{\left\vert{w}\right\vert
 ^{2}=-\rho (z)}{\left\vert{F(z,w)}\right\vert ^{p}d\left\vert{w}\right\vert
 }\rbrace dm(z)}=:{\left\Vert{F}\right\Vert}_{p}^{p}<\infty .$ \ \par 
but  $\displaystyle F(z,w)$  is holomorphic in  $w$  for  $\displaystyle
 z$  fixed, hence  $\displaystyle \ \left\vert{F(z,w)}\right\vert
 ^{p}$  is sub harmonic in  $w$  which implies\ \par 
\quad \quad \quad \quad   $\displaystyle \ \int_{\left\vert{w}\right\vert ^{2}=-\rho
 (z)}{\left\vert{F(z,w)}\right\vert ^{p}d\left\vert{w}\right\vert
 }\geq \left\vert{F(z,0)}\right\vert ^{p}.$ \ \par 
Hence\ \par 
\quad \quad \quad \quad \quad 	 $\displaystyle \ \int_{\Omega }{h(z)\left\vert{F(z,0)^{p}}\right\vert
 dm(z)}\leq {\left\Vert{F}\right\Vert}_{p}^{p}<\infty ,$ \ \par 
which implies, because  $\displaystyle h(z)$  is bounded below
 and above in  $\displaystyle \bar \Omega ,$  that\ \par 
\quad \quad \quad \quad \quad 	 $\displaystyle \ \int_{\Omega }{\left\vert{F(z,0)}\right\vert
 ^{p}dm(z)}\leq \frac{1}{\alpha }{\left\Vert{F}\right\Vert}_{p}^{p}<\infty
 .$ \ \par 
\quad \quad  	Now apply this for  $\tilde \Omega _{\epsilon }$  instead of
  $\tilde \Omega \ ;$  we have that  $\displaystyle F(z,w)$ 
 is continuous up to  $\displaystyle \partial \tilde \Omega _{\epsilon
 }$  because  $\epsilon >0.$  So\ \par 
\quad \quad \quad  $\displaystyle \ \int_{\partial \tilde \Omega _{\epsilon }}{\left\vert{F(z,w)}\right\vert
 ^{p}d\tilde \sigma (z,w)}\geq \alpha \int_{\Omega _{\epsilon
 }}{\left\vert{F(z,0)}\right\vert ^{p}dm(z)}.$ \ \par 
Hence by Fatou's lemma with  $\epsilon \rightarrow 0,$ \ \par 
\quad \quad \quad  $\displaystyle \alpha {\left\Vert{F(\cdot ,0)}\right\Vert}_{A^{p}(\Omega
 )}\leq {\left\Vert{F}\right\Vert}_{H^{p}(\tilde \Omega )}^{p}.$ \ \par 
So we have the first part of the lemma.\ \par 
\quad  Conversely if  $\displaystyle f\in A^{p}(\Omega ),$  setting
  $F(z,w):=f(z)$  and reversing the previous computations, using
 equalities this time,\ \par 
\quad \quad \quad \quad \quad 	 $\displaystyle \ \int_{\partial \tilde \Omega }{\left\vert{F}\right\vert
 ^{p}d\tilde \sigma }=\int_{\Omega }{h(z)\lbrace \int_{\left\vert{w}\right\vert
 ^{2}=-\rho (z)}{\left\vert{F(z,w)}\right\vert ^{p}d\left\vert{w}\right\vert
 }\rbrace dm(z)}=\int_{\Omega }{h(z)\left\vert{f(z)}\right\vert
 ^{p}dm(z)},$ \ \par 
because  $\displaystyle \ \int_{\left\vert{w}\right\vert ^{2}=-\rho
 (z)}{d\left\vert{w}\right\vert }=1.$  Hence\ \par 
\quad \quad \quad \quad \quad 	 $\displaystyle \ \int_{\partial \tilde \Omega }{\left\vert{F}\right\vert
 ^{p}d\tilde \sigma }\leq \beta \int_{\Omega }{\left\vert{f(z)}\right\vert
 ^{p}dm(z)}=\beta {\left\Vert{f}\right\Vert}_{A^{p}(\Omega )}^{p}.$
   $\blacksquare $ \ \par 
The only thing we used was that  $\displaystyle \ \left\vert{F(z,w)}\right\vert
 ^{p}$  is sub harmonic in  $w$  for  $z$  fixed. This being
 also true for   $\displaystyle F(z,w)\in {\mathcal{N}}(\tilde
 \Omega ),$  the very same proof gives\ \par 
\begin{Proposition}
 Let  $\displaystyle \Omega $  be a domain in  ${\mathbb{C}}^{n}$
  and  $\displaystyle \tilde \Omega $  its lift in  ${\mathbb{C}}^{n+1}.$
  There are constants  $\alpha >0,\ \beta >0$  depending only
 on  $\displaystyle \Omega $  such that if  $\displaystyle F\in
 {\mathcal{N}}(\tilde \Omega ),$  then  $F(z,0)\in {\mathcal{N}}_{0}(\Omega
 )$  and  $\ {\left\Vert{F(\cdot ,0)}\right\Vert}_{{\mathcal{N}}_{0}(\Omega
 )}\leq \frac{1}{\alpha }{\left\Vert{F}\right\Vert}_{{\mathcal{N}}(\tilde
 \Omega )}.$ \par 
Conversely if  $f\in {\mathcal{N}}_{0}(\Omega ),$  for  $\displaystyle
 (z,w)\in \tilde \Omega $  set  $\displaystyle F(z,w):=f(z),$
  then we have  $\ {\left\Vert{F}\right\Vert}_{{\mathcal{N}}(\tilde
 \Omega )}\leq \beta {\left\Vert{f}\right\Vert}_{{\mathcal{N}}_{0}(\Omega
 )}.$ \par 
\end{Proposition}
\quad \quad  	Now if we start with a function  $\displaystyle F(z,w)\in A^{p}(\tilde
 \Omega )$  what happens ? We have\ \par 
\begin{Proposition}
 ~\label{subPrinFin1455}Let  $\displaystyle \Omega $  be a domain
 in  ${\mathbb{C}}^{n}$  and  $\displaystyle \tilde \Omega $
  its lift in  ${\mathbb{C}}^{n+1}.$  If  $\displaystyle F\in
 A^{p}(\tilde \Omega ),$  then\par 
\quad \quad \quad \quad   $F(z,0)\in A^{p}_{1}(\Omega )$  and  $\ {\left\Vert{F(\cdot
 ,0)}\right\Vert}_{A_{1}^{p}(\Omega )}\leq \frac{1}{\pi }{\left\Vert{F}\right\Vert}_{A^{p}(\tilde
 \Omega )}.$ \par 
Conversely if  $f\in A_{1}^{p}(\Omega ),$  for  $\displaystyle
 (z,w)\in \tilde \Omega $  set  $\displaystyle F(z,w):=f(z),$
  then we have  $\ {\left\Vert{F}\right\Vert}_{A^{p}(\tilde \Omega
 )}\leq \pi {\left\Vert{f}\right\Vert}_{A_{1}^{p}(\Omega )}.$ \par 
\end{Proposition}
\quad \quad  	Proof.\ \par 
By Fubini we have\ \par 
\quad \quad \quad \quad \quad 	 $\displaystyle \ \int_{\tilde \Omega }{\left\vert{F(z,w)}\right\vert
 ^{p}dm(z,w)}=\int_{\Omega }{\lbrace \int_{\left\vert{w}\right\vert
 ^{2}<-r(z)}{\left\vert{F(z,w)}\right\vert ^{p}dm(w)}\rbrace dm(z)}.$ \ \par 
But again  $\displaystyle \ \left\vert{F(z,w)}\right\vert ^{p}$
  is sub harmonic in  $w$  for  $z$  fixed hence\ \par 
\quad \quad \quad \quad \quad 	 $\displaystyle \ \left\vert{F(z,0)}\right\vert ^{p}\leq \frac{1}{\pi
 (-r(z))}\int_{\left\vert{w}\right\vert ^{2}<-r(z)}{\left\vert{F(z,w)}\right\vert
 ^{p}dm(w)},$ \ \par 
because the area of the disc  $\displaystyle \ \lbrace \left\vert{w}\right\vert
 ^{2}<-r(z)\rbrace $  is  $\displaystyle \pi (-r(z)).$  So\ \par 
\quad \quad \quad \quad \quad 	 $\displaystyle \pi \int_{\Omega }{\left\vert{F(z,0)}\right\vert
 ^{p}(-r(z))dm(z)\leq }\int_{\tilde \Omega }{\left\vert{F(z,w)}\right\vert
 ^{p}dm(z,w)}$ \ \par 
hence\ \par 
\quad \quad \quad \quad \quad 	 $\displaystyle \ {\left\Vert{F(\cdot ,0)}\right\Vert}_{A_{1}^{p}(\Omega
 )}\leq \frac{1}{\pi }{\left\Vert{F}\right\Vert}_{A^{p}(\tilde
 \Omega )}.$ \ \par 
\quad \quad  	Conversely if  $\displaystyle F(z,w)=f(z)\in A_{1}^{p}(\Omega ),$ \ \par 
\quad \quad \quad \quad \quad 	 $\displaystyle \ \int_{\tilde \Omega }{\left\vert{F(z,w)}\right\vert
 ^{p}dm(z,w)}=\int_{\Omega }{\left\vert{f(z)}\right\vert ^{p}\lbrace
 \int_{\left\vert{w}\right\vert ^{2}<-r(z)}{dm(w)}\rbrace dm(z)}=\int_{\Omega
 }{\left\vert{f(z)}\right\vert ^{p}\pi (-r(z))dm(z)}$ \ \par 
hence\ \par 
\quad \quad \quad \quad \quad 	 $\displaystyle \ {\left\Vert{F}\right\Vert}_{A^{p}(\tilde \Omega
 )}\leq \pi {\left\Vert{f}\right\Vert}_{A_{1}^{p}(\Omega )}.$
   $\blacksquare $ \ \par 
\quad \quad  	We have the same results with the same proofs replacing Bergman
 classes by Nevanlinna ones.\ \par 
\begin{Proposition}
 Let  $\displaystyle \Omega $  be a domain in  ${\mathbb{C}}^{n}$
  and  $\displaystyle \tilde \Omega $  its lift in  ${\mathbb{C}}^{n+1}.$
  If  $F\in {\mathcal{N}}_{0}(\tilde \Omega ),$  then\par 
\quad \quad \quad \quad \quad   $F(z,0)\in {\mathcal{N}}_{1}(\Omega )$  and  $\ {\left\Vert{F(\cdot
 ,0)}\right\Vert}_{{\mathcal{N}}_{1}(\Omega )}\leq \frac{1}{\pi
 }{\left\Vert{F}\right\Vert}_{{\mathcal{N}}_{0}(\tilde \Omega )}.$ \par 
Conversely if  $f\in {\mathcal{N}}_{1}(\Omega ),$  for  $\displaystyle
 (z,w)\in \tilde \Omega $  set  $\displaystyle F(z,w):=f(z),$
  then we have  $\ {\left\Vert{F}\right\Vert}_{{\mathcal{N}}_{0}(\tilde
 \Omega )}\leq \pi {\left\Vert{f}\right\Vert}_{{\mathcal{N}}_{1}(\Omega
 )}.$ \par 
\end{Proposition}
\quad \quad  	Proof of the subordination lemma.\ \par 
We prove the subordination lemma for a one level lift. To get
 it for  $k$  levels lift, we just proceed by induction remarking that\ \par 
\quad \quad \quad \quad   $\displaystyle \widetilde{(\tilde \Omega _{k-1})}=\tilde \Omega
 _{k}.$ 	 $\displaystyle $ \ \par 
\quad  Let  $\displaystyle \Omega $  be a domain in  ${\mathbb{C}}^{n}$
  and set  $\displaystyle \tilde \Omega _{k}$  its  $k$  steps
 lift. Let  $\displaystyle F(z,w_{1},...,w_{k})\in H^{p}(\tilde
 \Omega _{k})$  then by the one level lift, proposition~\ref{subPrinFin1454}
 we have\ \par 
\quad \quad \quad \quad \quad 	 $\displaystyle F(z,w_{1},...,w_{k-1},0)\in A^{p}(\tilde \Omega
 _{k-1}),\ {\left\Vert{F(\cdot ,0)}\right\Vert}_{A^{p}(\tilde
 \Omega _{k-1})}\leq \frac{1}{\alpha }{\left\Vert{F}\right\Vert}_{H^{p}(\tilde
 \Omega _{k})}.$ \ \par 
Now set  $\displaystyle F_{1}(z,w_{1},...,w_{k-1}):=F(z,w_{1},...,w_{k-1},0)\in
 A^{p}(\tilde \Omega _{k-1})$  and apply proposition~\ref{subPrinFin1455},
 we get\ \par 
\quad \quad 	 $\displaystyle F_{1}(z,w_{1},...,w_{k-2},0)\in A_{1}^{p}(\tilde
 \Omega _{k-2}),\ {\left\Vert{F_{1}(\cdot ,0)}\right\Vert}_{A_{1}^{p}(\tilde
 \Omega _{k-2})}\leq \frac{1}{\pi }{\left\Vert{F_{1}}\right\Vert}_{A^{p}(\tilde
 \Omega _{k-1})}\leq \frac{1}{\alpha \pi }{\left\Vert{F}\right\Vert}_{H^{p}(\tilde
 \Omega _{k})}.$ \ \par 
And so on.\ \par 
The converse is done the same way as for the Nevanlinna classes.
  $\blacksquare $ \ \par 
\quad  Exactly the same induction gives the easy corollary :\ \par 
\begin{Corollary}
 Let  $\Omega $  be a domain in  ${\mathbb{C}}^{n},\ \tilde \Omega
 $  its lift in  ${\mathbb{C}}^{n+k}$  and  $F(z,w)\in A_{l}^{p}(\tilde
 \Omega ),$  we have  $f(z):=F(z,0)\in A_{k+l}^{p}(\Omega )$
  and  $\ {\left\Vert{f}\right\Vert}_{A_{k+l}^{p}(\Omega )}\lesssim
 {\left\Vert{F}\right\Vert}_{A_{l}^{p}(\tilde \Omega )};$ \par 
 if  $F(z,w)\in {\mathcal{N}}_{l}(\tilde \Omega ),$  then  $f(z):=F(z,0)\in
 {\mathcal{N}}_{k+l}(\Omega )$  and  $\ {\left\Vert{f}\right\Vert}_{{\mathcal{N}}_{k+l}(\Omega
 )}\lesssim {\left\Vert{F}\right\Vert}_{{\mathcal{N}}_{l}(\tilde
 \Omega )}.$ \par 
 A function  $f,$  holomorphic in  $\Omega ,$  is in the Bergman
 space  $A_{k+l}^{p}(\Omega )$  (resp. in the Nevanlinna Bergman
 space  ${\mathcal{N}}_{k+l}(\Omega )$ ) if and only if the function
  $F(z,w):=f(z)$  is in the Bergman space  $A_{l}^{p}(\tilde
 \Omega )$  (resp. in the Nevanlinna class  ${\mathcal{N}}_{l}(\tilde
 \Omega )$ ) and we have  $\ {\left\Vert{f}\right\Vert}_{A_{k+l}^{p}}\simeq
 {\left\Vert{F}\right\Vert}_{A_{l}^{p}(\tilde \Omega )}$  (resp.
  $\ {\left\Vert{f}\right\Vert}_{{\mathcal{N}}_{k+l}(\Omega )}\simeq
 {\left\Vert{F}\right\Vert}_{{\mathcal{N}}_{l}(\tilde \Omega )}$ ).\par 
\end{Corollary}

\section{Geometric Carleson measures and  $p$ -Carleson measures.~\label{subPrinFin435}}

      In order to define precisely the geometric Carleson measures,
 we need the notion of a "good" family of polydiscs, directly
  inspired by the work of Catlin~\cite{Catlin84} and introduced
 in~\cite{AmAspc09}.\ \par 
\quad  Let  ${\mathcal{U}}$  be a neighbourhood of  $\partial \Omega
 $  in  $\Omega $  such that the normal projection  $\pi $  onto
  $\partial \Omega $  is a smooth well defined application.\ \par 
\quad  Let  $\alpha \in \partial \Omega $  and let  $b(\alpha )=(L_{1},\
 L_{2},...,\ L_{n})$  be an orthonormal basis of  ${\mathbb{C}}^{n}$
  such that  $(L_{2},...,\ L_{n})$  is a basis of the tangent
 complex space  $T_{\alpha }^{{\mathbb{C}}}$  of  $\partial \Omega
 $  at  $\alpha \ ;$  hence  $L_{1}$  is the complex normal at
  $\alpha $  to  $\partial \Omega .$ \ \par 
Let  $m(\alpha )=(m_{1},\ m_{2},...,\ m_{n})\in {\mathbb{N}}^{n}$
  be a multi-index at  $\alpha $  with  $m_{1}=1,\ \forall j\geq
 2,\ m_{j}\geq 2.$ \ \par 
For  $a\in {\mathcal{U}}$  and  $t>0$  set  $\alpha =\pi (a)$
  and  $P_{a}(t):=\prod_{j=1}^{n}{tD_{j}},$  the polydisc such
 that  $tD_{j}$  is the disc centered at  $a,$  parallel to 
 $L_{j}\in b(\alpha ),$  with radius  $t\left\vert{r(a)}\right\vert
 ^{1/m_{j}}$  (recall that we have  $\ \left\vert{r(a)}\right\vert
 \simeq \delta (a)$  ).\ \par 
\quad  Set  $b(a):=b(\pi (a)),\ m(a):=m(\pi (a)),$  for  $\displaystyle
 a\in {\mathcal{U}}.$ \ \par 
\quad  This way we have a family of polydiscs  ${\mathcal{P}}:=\lbrace
 P_{a}(t)\rbrace _{a\in {\mathcal{U}}}$  defined by the family
 of basis  $\lbrace b(a)\rbrace _{a\in {\mathcal{U}}},$  the
 family of multi-indices  $\lbrace m(a)\rbrace _{a\in {\mathcal{U}}}$
  and the number  $t.$  Notice that the polydisc  $P_{a}(2)$
  always overflows the domain  $\Omega .$ \ \par 
\quad  It will be useful to extend this family to the whole of  $\Omega
 .$  In order to do so let  $(z_{1},\ ...,\ z_{n})$  be the canonical
 coordinates system in  ${\mathbb{C}}^{n}$  and for  $a\in \Omega
 \backslash {\mathcal{U}},$  let  $P_{a}(t)$  be the polydisc
 of center  $a,$  of sides parallel to the axis and radius  $t\delta
 (a)$  in the  $z_{1}$  direction and  $t\delta (a)^{1/2}$  in
 the other directions. So the points  $a\in \Omega \backslash
 {\mathcal{U}}$  have automatically a "minimal" multi-index 
 $m(a)=(1,\ 2,\ ...,\ 2).$ \ \par 
\quad  Now we can set\ \par 
\begin{Definition}
 We say that  ${\mathcal{P}}$  is a "good family" of polydiscs
 for  $\Omega $  if the  $m_{j}(a)$  are uniformly bounded on
  $\Omega $  and if it exists  $\delta _{0}>0$  such that all
 the polydiscs  $P_{a}(\delta _{0})$  of  ${\mathcal{P}}$  are
 contained in  $\Omega .$  In this case we call  $m(a)$  the
 multi-type at  $a$  of the family  ${\mathcal{P}}.$ \par 
\end{Definition}
\quad  We notice that, for a good family  ${\mathcal{P}},$  by definition
 the multi-type is always finite. Moreover there is no regularity
 assumptions on the way that the basis  $b(a)$  varies with respect
 to  $a\in \Omega .$ \ \par 
\quad  We can see easily that there are always good families of polydiscs
 in a domain  $\Omega $  in  ${\mathbb{C}}^{n}\ :$  for a point
  $a\in \Omega ,$  take any orthonormal basis  $b(a)=(L_{1},\
 L_{2},...,\ L_{n}),$  with  $L_{1}$  the complex normal direction,
 and the "minimal" multitype  $m(a)=(1,\ 2,...,\ 2).$  Then,
 because the level sets  $\partial \Omega _{a}$  are uniformly
 of class  ${\mathcal{C}}^{2}$  and compact, we have the existence
 of a uniform  $\delta _{0}>0$  such that the family  ${\mathcal{P}}$
  is a good one. As seen in~\cite{AmAspc09}, in the  strictly
 pseudo-convex domains, this family with "minimal" multi-type
 is the right one.\ \par 
\quad  We can give the definitions relative to Carleson measures.\ \par 
\begin{Definition}
 A positive borelian  measure  $\mu $  on  $\Omega $  is a geometric
 Carleson measure,  $\mu \in \Lambda (\Omega ),$  if\par 
\quad \quad \quad  $\displaystyle \exists C=C_{\mu }>0::\forall a\in \Omega ,\
 \mu (\Omega \cap P_{a}(2))\leq C\sigma (\partial \Omega \cap P_{a}(2)).$ \par 
\end{Definition}
\begin{Definition}
 A positive borelian measure  $\mu $  on  $\Omega $  is a  $p$
 -Carleson measure in  $\Omega $  if\par 
\quad \quad \quad  $\displaystyle \exists C>0::\forall f\in H^{p}(\Omega ),\ \int_{\Omega
 }{\left\vert{f(z)}\right\vert ^{p}\,d\mu (z)}\leq C^{p}{\left\Vert{f}\right\Vert}_{H^{p}(\Omega
 )}^{p}.$ \par 
\end{Definition}
\quad  And  analogously for the Bergman spaces.\ \par 
\begin{Definition}
 A positive borelian  measure  $\mu $  on  $\Omega $  is a  $k$
 -geometric Bergman-Carleson measure,  $\mu \in \Lambda _{k}^{}(\Omega
 ),$  if\par 
\quad \quad \quad  $\displaystyle \exists C=C_{\mu }>0::\forall a\in \Omega ,\
 \mu (\Omega \cap P_{a}(2))\leq Cm_{k-1}(\Omega \cap P_{a}(2)).$ \par 
\end{Definition}
Notice the gap  $k\rightarrow k-1.$ \ \par 
\begin{Definition}
 A positive borelian measure  $\mu $  is  $(p,k)$ -Bergman-Carleson
 measure in  $\Omega $  if\par 
\quad \quad \quad  $\displaystyle \exists C>0::\forall f\in A_{k-1}^{p}(\Omega
 ),\ \int_{\Omega }{\left\vert{f(z)}\right\vert ^{p}\,d\mu (z)}\leq
 C^{p}{\left\Vert{f}\right\Vert}_{A_{k-1}^{p}(\Omega )}^{p}.$ \par 
\end{Definition}
\ \par 
\begin{Definition}
 We shall say that the domain  $\Omega $  has the  $p$ -Carleson
 embedding property,  $p$ -CEP, if\par 
\quad \quad \quad  $\displaystyle \forall \mu \in \Lambda (\Omega ),\ \exists C=C_{\mu
 }>0::\forall f\in H^{p}(\Omega ),\ \int_{\Omega }{\left\vert{f}\right\vert
 ^{p}\,d\mu }\leq C{\left\Vert{f}\right\Vert}_{H^{p}(\Omega )}^{p}.$ \par 
\end{Definition}
\quad  And the same for the Bergman spaces.\ \par 
\begin{Definition}
 We shall say that the domain  $\Omega $  has the  $(p,k)$ -Bergman-Carleson
 embedding property,  $(p,\ k)$ -BCEP, if\par 
\quad \quad \quad  $\displaystyle \forall \mu \in \Lambda _{k}(\Omega ),\ \exists
 C=C_{\mu ,p}>0::\forall f\in A_{k-1}^{p}(\Omega ),\ \int_{\Omega
 }{\left\vert{f}\right\vert ^{p}\,d\mu }\leq C{\left\Vert{f}\right\Vert}_{A_{k-1}^{p}(\Omega
 )}^{p}.$ \par 
\end{Definition}

\subsection{The subordination lemma applied to Carleson measures.~\label{subPrinFin436}}
\
 \par 
\quad  We shall fix  $k\in {\mathbb{N}}$  and lift the  measure on
 the domain  $\tilde \Omega :=\lbrace \tilde r(z,w):=r(z)+\left\vert{w}\right\vert
 ^{2}<0\rbrace ,$  with  $w=(w_{1},...,w_{k})\in {\mathbb{C}}^{k}.$
  We already know how to lift a function, the lifted measure
  $\tilde \mu $  of a measure  $\mu $  is just\ \par 
\quad \quad \quad  $\tilde \mu :=\mu \otimes \delta ,$ \ \par 
with  $\delta $  the delta Dirac measure of the origin in  ${\mathbb{C}}^{k}.$
  We shall need a lemma linking Bergman and Hardy geometric Carleson
 measures.\ \par 
\quad  Let  $\Omega $  be a domain in  ${\mathbb{C}}^{n},\ \tilde \Omega
 $  be its lift in  ${\mathbb{C}}^{n+k},$  and suppose that 
 $\tilde \Omega $  is equipped with a good family of polydiscs
  $\tilde {\mathcal{P}},$  we have the definition :\ \par 
\begin{Definition}
 We shall say that the good family of polydiscs  $\tilde {\mathcal{P}}$
  on the domain  $\tilde \Omega $  is "homogeneous" if\par 
\quad \quad  	(Hg)       $\exists t>0,\ \exists C>0::\forall a\in \tilde
 \Omega ,\ \Omega \cap \tilde P_{a}(2)\neq \emptyset ,$ \par 
\quad \quad \quad \quad \quad 	 $\displaystyle \ \forall b\in \Omega \cap \tilde P_{a}(2),\
 \tilde P_{b}(t)\supset \tilde P_{a}(2)\ and\ \tilde \sigma (\partial
 \tilde \Omega \cap \tilde P_{b}(t))\leq C\tilde \sigma (\partial
 \tilde \Omega \cap \tilde P_{a}(2)),$ \par 
where  $\Omega =\tilde \Omega \cap \lbrace w=0\rbrace \subset
 \tilde \Omega .$ \par 
\end{Definition}
\quad  Naturally the domain  $\Omega $  is equipped with the family
  ${\mathcal{P}}$  induced by  $\tilde {\mathcal{P}}$  the following way\ \par 
\quad \quad \quad  $\forall a\in \Omega ,\ P_{a}(u):=\tilde P_{(a,0)}(u)\cap \lbrace
 w=0\rbrace ,$ \ \par 
which is easily seen to be a good family for  $\Omega .$ \ \par 
As examples we have the  strictly pseudo-convex domains and
 the convex domains of finite type, because both are domains
 of homogeneous type in the sense of Coifman-Weiss~\cite{CoifWeiss71}.\ \par 
\begin{Lemma}
 ~\label{subPrinciple621}Let  $(\Omega ,\ \tilde \Omega )$  be
 as above and suppose that  $\tilde \Omega $  is equipped with
 a good family of polydiscs  $\tilde {\mathcal{P}}$  which verifies
 the hypothesis (Hg). The measure  $\mu $  is a  $k$ -geometric
 Bergman-Carleson measure in  $\Omega $  iff the measure  $\tilde
 \mu $  is a geometric Carleson measure in  $\tilde \Omega .$ \par 
\end{Lemma}
\quad  Proof.\ \par 
Suppose that  $\mu $  is a  $k$ -geometric Bergman-Carleson
 measure in  $\Omega ,$  we want to show :\ \par 
\quad \quad \quad  $\displaystyle \exists C>0::\forall (a,\ b)\in \tilde \Omega
 ,\ \tilde \mu (\tilde \Omega \cap \tilde P_{(a,b)}(2))\leq C\tilde
 \sigma (\partial \tilde \Omega \cap \tilde P_{(a,b)}(2)),$ \ \par 
with  $\tilde P_{c}$  the polydisc of center  $c=(a,b)\in \tilde
 \Omega $  of the family  $\tilde {\mathcal{P}}.$  Let us see
 first the case where  $b=0,$  i.e.  $(a,b)=(a,0)\in \Omega \subset
 \tilde \Omega .$  Then, by definition of  $\tilde \mu ,$  we have\ \par 
\quad \quad \quad  $\displaystyle \tilde \mu (\tilde \Omega \cap \tilde P_{(a,0)}(2))=\mu
 (\Omega \cap P_{a}(2)).$ \ \par 
On the other hand, we have, exactly as in the proof of the subordination
 lemma,\ \par 
\quad \quad \quad  $\displaystyle \tilde \sigma (\partial \tilde \Omega \cap \tilde
 P_{(a,0)}(2))\simeq \int_{\Omega \cap P_{a}(2)}{kv_{k}(-r(z))^{k-1}\,dm(z)}=m_{k-1}(\Omega
 \cap P_{a}(2)).$ \ \par 
But if  $\mu $  is a  $k$ -geometric Bergman-Carleson measure
 in  $\Omega ,$  we have\ \par 
\quad \quad \quad  $\displaystyle \exists C>0::\forall a\in \Omega ,\ \mu (\Omega
 \cap P_{a}(2))\leq Cm_{k-1}(\Omega \cap P_{a}(2)),$ \ \par 
so\ \par 
\quad \quad \quad  $\displaystyle \tilde \mu (\tilde \Omega \cap \tilde P_{(a,0)}(2))=\mu
 (\Omega \cap P_{a}(2))\leq Cm_{k-1}(\Omega \cap P_{a}(2))\simeq
 C\tilde \sigma (\partial \tilde \Omega \cap \tilde P_{(a,0)}(2)).$ \ \par 
Now take a general  $\tilde P_{(a,b)}(2).$  In order for  $\tilde
 \mu (\tilde \Omega \cap \tilde P_{(a,b)}(2))$  to be non zero,
 we must have\ \par 
\quad \quad \quad  $\tilde P_{(a,b)}(2)\cap \lbrace w=0\rbrace \neq \emptyset \Rightarrow
 \exists (c,0)\in \tilde P_{(a,b)}(2)\cap \lbrace w=0\rbrace .$ \ \par 
By the (Hg) hypothesis, this means that we have  $\tilde P_{(c,0)}(t)\supset
 \tilde P_{(a,b)}(2)$  with the uniform control  $\tilde \sigma
 (\partial \tilde \Omega \cap \tilde P_{(c,0)}(t))\lesssim \tilde
 \sigma (\partial \tilde \Omega \cap \tilde P_{(a,b)}(2)).$ \ \par 
We apply the above inequality\ \par 
\quad \quad \quad  $\displaystyle \tilde \mu (\tilde \Omega \cap \tilde P_{(a,b)}(2))\leq
 \tilde \mu (\tilde \Omega \cap \tilde P_{(c,0)}(t))\leq Cm_{k}(\Omega
 \cap P_{c}(t))=C\tilde \sigma (\partial \tilde \Omega \cap \tilde
 P_{(c,0)}(t))\lesssim \tilde \sigma (\partial \tilde \Omega
 \cap \tilde P_{(a,b)}(2)),$ \ \par 
hence  $\tilde \mu $  is a geometric Carleson measure on  $\tilde
 \Omega .$ \ \par 
\quad  Conversely suppose that  $\tilde \mu $  is a geometric Carleson
 measure on  $\tilde \Omega ,$  this means\ \par 
\quad \quad \quad  $\displaystyle \forall (a,b)\in \tilde \Omega ,\ \tilde \mu
 (\tilde \Omega \cap \tilde P_{(a,b)}(2))\leq C\tilde \sigma
 (\partial \tilde \Omega \cap \tilde P_{(a,b)}(2)),$ \ \par 
hence, in particular for  $b=0,$ \ \par 
\quad \quad \quad  $\forall a\in \Omega ,\ \tilde \mu (\tilde \Omega \cap P_{(a,0)}(2))\leq
 C\tilde \sigma (\tilde \Omega \cap P_{(a,0)}(2)),$ \ \par 
but then, by definition of  $\tilde \mu $  and with the previous
 computation of  $\displaystyle \tilde \sigma (\tilde \Omega
 \cap P_{(a,0)}(2)),$  we get\ \par 
\quad \quad \quad  $\displaystyle \forall a\in \Omega ,\ \mu (\Omega \cap P_{a}(2))\leq
 Cm_{k-1}(\Omega \cap P_{a}(2)),$ \ \par 
hence the measure  $\mu $  is a  $k$ -geometric Bergman-Carleson
 measure in  $\Omega .$   $\blacksquare $ \ \par 
\quad  Now we shall use the subordination lemma to get a Bergman-Carleson
 embedding theorem from a Hardy-Carleson embedding one.\ \par 
\begin{Theorem}
 ~\label{subPrinFin647}Let  $(\Omega ,\ \tilde \Omega )$  be
 as usual and suppose that  $\tilde \Omega $  is equipped with
 a good family of polydiscs  $\tilde {\mathcal{P}}$  which verifies
 the hypotheses (Hg). If the lifted domain  $\tilde \Omega $
  has the  $p$ -CEP then  $\Omega $  has the  $(p,k)$ -BCEP.\par 
\end{Theorem}
\quad  Proof.\ \par 
Suppose the positive measure  $\mu $  is a  $k$ -geometric Bergman-Carleson
 measure ; by the previous lemma, we have that the lifted measure
  $\tilde \mu $  is a geometric Carleson measure in  $\tilde
 \Omega .$  By the  $p$ -CEP we have\ \par 
\quad \quad \quad  $\forall F\in H^{p}(\tilde \Omega ),\ \int_{\tilde \Omega }{\left\vert{F}\right\vert
 ^{p}\,d\tilde \mu }\leq C_{\mu }^{p}{\left\Vert{F}\right\Vert}_{H^{p}(\tilde
 \Omega )}^{p}.$ \ \par 
Choose  $f(z)\in A_{k-1}^{p}(\Omega )$  and set  $\forall (z,w)\in
 \tilde \Omega ,\ F(z,w)=f(z).$  By the subordination lemma we have\ \par 
\quad \quad \quad \quad  $\ {\left\Vert{f}\right\Vert}_{A_{k-1}^{p}(\Omega )}\simeq {\left\Vert{F}\right\Vert}_{H^{p}(\tilde
 \Omega )},$ \ \par 
and by definition of  $\tilde \mu ,$  we have\ \par 
\quad \quad \quad  $\displaystyle \ \int_{\Omega }{\left\vert{f}\right\vert ^{p}\,d\mu
 }=\int_{\tilde \Omega }{\left\vert{F}\right\vert ^{p}\,d\tilde
 \mu }\leq C_{\mu }^{p}{\left\Vert{F}\right\Vert}_{H^{p}(\tilde
 \Omega )}^{p}\lesssim {\left\Vert{f}\right\Vert}_{A_{k-1}^{p}(\Omega
 )},$ \ \par 
hence  $\mu $  is a  $(k,p)$ -Bergman-Carleson measure in  $\Omega
 .$   $\blacksquare $ \ \par 
\begin{Theorem}
 ~\label{subPrinFin440}Let  $(\Omega ,\ \tilde \Omega )$  be
 as usual and suppose that  $\tilde \Omega $  is equipped with
 a good family of polydiscs  $\tilde {\mathcal{P}}$  which verifies
 the hypotheses (Hg). If  $p$ -Carleson implies geometric Carleson
 in  $\tilde \Omega ,$  then  $(p,k)$ -Bergman-Carleson implies
 geometric  $k$ -Bergman-Carleson in  $\Omega .$ \par 
\end{Theorem}
\quad  Proof.\ \par 
If the positive measure  $\mu $  is  $(p,k)$ -Bergman-Carleson
 in  $\Omega $  then  $\tilde \mu $  is a  $p$ -Carleson measure
 in  $\tilde \Omega $  by lemma~\ref{subPrinciple621} hence a
 geometric Carleson measure in  $\tilde \Omega $  by the assumption
 of the theorem. Then applying lemma~\ref{subPrinciple621} we
 get that  $\mu $  is a  $k$ -geometric Carleson measure in 
 $\Omega $  hence the theorem.  $\blacksquare $ \ \par 
\begin{Remark}
The definition of geometric Carleson measures depends on the
 chosen good family of polydiscs on the domain ; the theorem
 asserts the equivalence of properties between a domain  $\Omega
 $  and its lift  $\tilde \Omega .$  The fact that a lifted domain
  $\tilde \Omega $  equipped with a good family of polydiscs
  $\tilde {\mathcal{P}}$  has the Carleson embedding property
 has to be proved directly but if it has the  $p$ -CEP  then
  $\Omega $  equipped with the induced family  ${\mathcal{P}}$
  has the  $(p,k)$ -BCEP without any further proof.\par 
\end{Remark}

\subsection{Application to strictly pseudo-convex domains.}
\ \par 
\begin{Corollary}
Let  $\Omega $  be a strictly pseudo-convex domain equipped
 with its minimal good family of polydiscs, then  $\Omega $ 
 has the  $(p,k)$  Bergman Carleson embedding property.\par 
\end{Corollary}
\quad \quad  	Proof.\ \par 
The domain  $\Omega $  equipped with its minimal good family
 has the  $p$ -CEP by H\"ormander~\cite{HormPSH67}, hence we
 can apply theorem~\ref{subPrinFin440}.  $\blacksquare $ \ \par 
\quad \quad  	This corollary gives a characterization of the  $(p,k)$ -Bergman-Carleson
 measures of the strictly pseudo-convex domains. Let  $\Omega
 $  be a strictly pseudo-convex domain and  $\tilde \Omega $
  its lift in  ${\mathbb{C}}^{n+k}.$  Let  $\tilde {\mathcal{P}}$
  be its minimal good family of polydiscs in  $\tilde \Omega
 \ ;$  one can see easily that the induced family of polydiscs
  ${\mathcal{P}}$  on  $\Omega $  is again the minimal good family
 of polydiscs. Recall that  $\displaystyle \forall a\in \Omega
 ,\ \delta (a)=d(a,\partial \Omega )\ ;$  we have this characterization :\ \par 
\begin{Corollary}
A positive Borel measure  $\mu $  in a strictly pseudo-convex
 domain in  ${\mathbb{C}}^{n}$  is a  $(p,k)$ -Bergman-Carleson
 measure iff :\par 
\quad \quad \quad \quad \quad 	 $\forall a\in \Omega ,\ \mu (P_{a}(2))\lesssim \delta (a)^{n+k}.$ \par 
This means that it is a characterization of the measures such that\par 
\quad \quad \quad \quad \quad 	 $\ \forall p\geq 1,\ \forall f\in A_{k-1}^{p}(\Omega ),\ \int_{\Omega
 }{\left\vert{f}\right\vert ^{p}\,d\mu }\lesssim {\left\Vert{f}\right\Vert}_{A_{k-1}^{p}(\Omega
 )}.$ \par 
In particular this characterization is independent of  $p\geq 1.$ \par 
\end{Corollary}
\quad \quad  	Proof.\ \par 
Let  $\tilde \Omega $  be the lift of  $\Omega $  in  ${\mathbb{C}}^{n+k}$
  and  $\tilde \mu $  be the lift of  $\mu $  on  $\tilde \Omega .$ \ \par 
\quad \quad  	Suppose that  $\mu $  is a  $(p,k)$ -Bergman Carleson measure
 in  $\Omega ,$  then  $\tilde \mu $  is a  $p$ -Carleson measure
 in  $\tilde \Omega $  by lemma~\ref{subPrinciple621} then by
 a theorem of H\"ormander~\cite{HormPSH67} the  $p$ -Carleson
 measures are precisely the geometric ones in  $\tilde \Omega
 ,$  hence we have\ \par 
\quad \quad \quad \quad \quad 	 $\forall \tilde a\in \tilde \Omega ,\ \tilde \mu (\tilde \Omega
 \cap \tilde P_{\tilde a}(2))\lesssim \tilde \sigma (\partial
 \tilde \Omega \cap \tilde P_{\tilde a}(2)).$ \ \par 
Now let  $a\in \Omega ,\ \tilde a:=(a,0)\in \tilde \Omega $
  then a classical computation gives  $\tilde \sigma (\partial
 \tilde \Omega \cap \tilde P_{\tilde a}(2))\lesssim \tilde \delta
 (\tilde a)^{n+k}=\delta (a)^{n+k}.$  By the definition of  $\tilde
 \mu $  we have\ \par 
\quad \quad \quad \quad \quad 	 $\delta (a)^{n+k}\gtrsim \tilde \mu (\tilde \Omega \cap \tilde
 P_{\tilde a}(2))=\mu (\tilde P_{\tilde a}(2)\cap \Omega )=\mu
 (P_{a}(2)\cap \Omega ),$ \ \par 
so  $\forall a\in \Omega ,\ \mu (P_{a}(2)\cap \Omega )\lesssim
 \delta (a)^{n+k}.$ \ \par 
\quad \quad  	Now suppose that  $\forall a\in \Omega ,\ \mu (P_{a}(2)\cap
 \Omega )\lesssim \delta (a)^{n+k}$  then we have, by the definition
 of  $\tilde \mu ,$  with  $\displaystyle \tilde a:=(a,0)\in
 \tilde \Omega ,$ \ \par 
\quad \quad \quad \quad \quad 	 $\displaystyle \ \tilde \mu (\tilde \Omega \cap \tilde P_{\tilde
 a}(2))\leq \tilde \delta (a)^{n+k}\simeq \tilde \sigma (\partial
 \tilde \Omega \cap \tilde P_{\tilde a}(2)).$ \ \par 
Doing exactly as in the proof of lemma~\ref{subPrinciple621}
 we have the same inequality with a bigger constant for all 
 $\tilde a\in \tilde \Omega ,$  hence  $\tilde \mu $  is a geometric
 Carleson measure in  $\tilde \Omega .$  So by H\"ormander~\cite{HormPSH67},
  $\tilde \mu $  is a  $p$ -Carleson measure in  $\tilde \Omega
 $  hence we have the embedding\ \par 
\quad \quad \quad \quad \quad 	 $\forall F\in H^{p}(\tilde \Omega ),\ \int_{\tilde \Omega }{\left\vert{F}\right\vert
 ^{p}d\tilde \mu }\lesssim {\left\Vert{F}\right\Vert}_{H^{p}(\tilde
 \Omega )}.$ \ \par 
Now we take  $f\in A_{k-1}^{p}(\Omega )$  and we set  $\forall
 (z,w)\in \tilde \Omega ,\ F(z,w):=f(z)$  by the subordination
 lemma we have  $\ {\left\Vert{F}\right\Vert}_{H^{p}(\tilde \Omega
 )}\simeq {\left\Vert{f}\right\Vert}_{A_{k-1}^{p}(\Omega )}$  and\ \par 
\quad \quad \quad \quad \quad   $\ \int_{\tilde \Omega }{\left\vert{F}\right\vert ^{p}d\tilde
 \mu }=\int_{\Omega }{\left\vert{f}\right\vert ^{p}dm_{k-1}}\lesssim
 {\left\Vert{F}\right\Vert}_{H^{p}(\tilde \Omega )}^{p}\simeq
 {\left\Vert{f}\right\Vert}_{A_{k-1}^{p}(\Omega )}^{p}.$   $\blacksquare
 $ \ \par 
\quad \quad 	 Cima and Mercer~\cite{CimaMercer95} characterized the Carleson
 measures for the spaces  $A_{\alpha }^{p}(\Omega )$  for  $\Omega
 $  strictly pseudo-convex, and with  $\alpha \geq 0.$  In the
 case where  $\alpha $  is an integer we recover their characterization,
 because one has easily, when  $\Omega $  is a strictly pseudo-convex
 domain, that  $P_{a}(2)\cap \Omega \simeq W(\pi (a),\ \delta
 (a))$  where  $W(\zeta ,h)$  is the classical Carleson window
 in  $\Omega .$ \ \par 
\begin{Remark}
 In the case of the unit ball  $\Omega $  of  ${\mathbb{C}}^{n},\
 \tilde \Omega \subset {\mathbb{C}}^{n+1}$  N. Varopoulos indicated
 me an alternative proof for the fact that  $F(z,w)\in H^{p}(\tilde
 \Omega )\Rightarrow F(z,0)\in A_{p}(\Omega )\ :$  the Lebesgue
 measure on  $\lbrace w=0\rbrace \cap \tilde \Omega $  is easily
 seen to be a geometric Carleson measure in  $\tilde \Omega ,$
  hence by the Carleson-H\"ormander embedding theorem~\cite{HormPSH67}
 we have\par 
\quad \quad \quad \quad \quad  $\displaystyle \ \int_{\Omega }{\left\vert{F(z,0)}\right\vert
 ^{p}\,dm(z)}\leq C{\left\Vert{F}\right\Vert}_{H^{p}(\tilde \Omega )},$ \par 
and the assertion. Of course this is still valid in codimension
  $k\geq 1,$  with the weighted Lebesgue measure on  $\Omega
 ,$  and for strictly pseudo-convex domains because the Carleson-H\"ormander
 embedding theorem is still valid there.  But this is just one
 direction of the lemma, it works only if there is a Carleson
 embedding theorem and this proof is much less elementary than
 the previous one.\par 
\quad  In fact we can reverse things and say that one part of the subordination
 lemma asserts that the weighted Lebesgue measure on  $\Omega
 $  is always a Carleson measure in  $\tilde \Omega ,\ \Omega
 $   strictly pseudo-convex or not.\par 
\end{Remark}

\subsection{Application to convex domains of finite type in
  ${\mathbb{C}}^{n}$ }
\ \par 
\quad  In~\cite{AmAspc09} we prove a Carleson embedding theorem for
 the convex domains of finite type in  ${\mathbb{C}}^{n}.$ \ \par 
\begin{Theorem}
 Let  $\Omega $  be a  convex domain of finite type in  ${\mathbb{C}}^{n}$
  ; if the measure  $\mu $  is a geometric Carleson measure we have\par 
\quad \quad \quad  $\forall p>1,\ \exists C_{p}>0,\ \forall f\in H^{p}(\Omega ),\
 \int_{\Omega }{\left\vert{f}\right\vert ^{p}\,d\mu }\leq C_{p}^{p}{\left\Vert{f}\right\Vert}_{H^{p}}^{p}.$
 \par 
\quad  Conversely if the positive measure  $\mu $  is  $p$ -Carleson
 for a  $p\in \lbrack 1,\ \infty \lbrack ,$  then it is a geometric
 Carleson measure, hence it is  $q$ -Carleson for any  $q\in
 \rbrack 1,\ \infty \lbrack .$ \par 
\end{Theorem}
\quad  We already know that if  $\Omega $  is a convex domain of finite
 type, so is  $\tilde \Omega $  with the same type. Moreover
 the hypothesis (Hg) is true for these domains equipped with
 a (slightly modified) McNeal family of polydiscs, so we can
 apply what precedes in this case to get from the Carleson embedding
 theorem the Bergman-Carleson embedding one.\ \par 
\begin{Theorem}
 Let  $\Omega $  be a  convex domain of finite type in  ${\mathbb{C}}^{n}$
  ; if the measure  $\mu $  is a  $k$ -geometric Bergman-Carleson
 measure, i.e.\par 
\quad \quad \quad \quad \quad 	 $\displaystyle \exists C>0::\forall a\in \Omega ,\ \mu (\Omega
 \cap P_{a}(2))\leq Cm_{k-1}(\Omega \cap P_{a}(2)),$ \par 
we have\par 
\quad \quad \quad  $\forall p>1,\ \exists C_{p}>0,\ \forall f\in A_{k-1}^{p}(\Omega
 ),\ \int_{\Omega }{\left\vert{f}\right\vert ^{p}\,d\mu }\leq
 C_{p}^{p}{\left\Vert{f}\right\Vert}_{A_{k-1}^{p}(\Omega )}^{p}.$ \par 
\quad  Conversely if the positive measure  $\mu $  is  $(p,k)$ -Bergman-Carleson
 for a  $p\in \lbrack 1,\ \infty \lbrack ,$  then it is a  $k$
 -geometric Bergman-Carleson measure, hence it is  $(q,k)$ -Bergman-Carleson
 for any  $q\in \rbrack 1,\ \infty \lbrack .$ \par 
\end{Theorem}

\section{Interpolating sequences for Bergman spaces.~\label{subPrinFin438}}

\subsection{On Bergman and Szeg\"o projections.}
\ \par 
\quad  Let  $\Omega $  be a domain in  ${\mathbb{C}}^{n},$  recall
 the definition of its Szeg\"o projection : this is the orthogonal
 projection  $P$  from  $L^{2}(\partial \Omega )$  onto  $H^{2}(\Omega
 )\ ;$  we shall note its kernel by  $S(z,\ \zeta ),$  i.e.\ \par 
\quad \quad \quad \quad \quad   $\displaystyle \forall f\in L^{2}(\partial \Omega ),\ Pf(z)=\int_{\partial
 \Omega }{S(z,w)f(\zeta )d\sigma (\zeta )}.$ \ \par 
\quad \quad  	The same way, recall the definition of the Bergman projection
 : this is the orthogonal projection  $P_{k}$  from  $L^{2}(\Omega
 ,\ dm_{k})$   onto  $A_{k}^{2}(\Omega ),$  the holomorphic functions
 on  $\Omega $  still in  $L^{2}(\Omega ,\ dm_{k}).$  We shall
 note its kernel by  $B_{k}(z,\ \zeta )$  i.e.\ \par 
\quad \quad \quad \quad \quad   $\displaystyle \forall f\in L^{2}(\Omega ,\ dm_{k}),\ P_{k}f(z)=\int_{\Omega
 }{B_{k}(z,w)f(\zeta )dm_{k}(\zeta )}.$ \ \par 
Let  $\tilde \Omega $  be the lifted domain of  $\Omega $  in
  ${\mathbb{C}}^{n+k}$  ; we shall use the notation\ \par 
\quad \quad \quad \quad \quad 	 $\forall z\in \Omega ,\ \tilde z:=(z,0)\in \tilde \Omega .$ \ \par 
\begin{Corollary}
 ~\label{subPrinCarl23}For any  $a\in \Omega $  the Bergman kernel
  $B_{k-1}(z,a)$  and the Szeg\"o kernel  $\tilde S((z,w),\ \tilde
 a)$  for the lifted domain  $\tilde \Omega ,$  verify,\par 
\quad \quad \quad  $\forall a\in \Omega ,\ \forall z\in \Omega ,\ B_{k-1}(z,a)=\tilde
 S(\tilde z,{\tilde a}).$ \par 
Moreover we have\par 
\quad \quad \quad  $\forall a\in \Omega ,\ {\left\Vert{B_{k-1}(\cdot ,\ a)}\right\Vert}_{A_{k-1}^{p}(\Omega
 )}\simeq {\left\Vert{\tilde S(\cdot ,\ \tilde a)}\right\Vert}_{H^{p}(\tilde
 \Omega )}.$ \par 
\end{Corollary}
\quad  Proof.\ \par 
Let  $f\in A(\Omega )$  be a holomorphic function in  $\Omega
 ,$  continuous up to  $\partial \Omega .$  Let\ \par 
\quad \quad \quad \quad  $\forall (z,w)\in \tilde \Omega ,\ F(z,w):=f(z).$ \ \par 
We have\ \par 
\quad \quad \quad  $\ \int_{\Omega }{f(z)\bar B_{k-1}(z,a)\,dm_{k-1}(z)}=f(a)=F(a,0)=\int_{\partial
 \tilde \Omega }{F(z,w){\overline{\tilde S}}((z,w),\ {\tilde
 a})\,d\sigma (z,w)},$ \ \par 
by the reproducing property of these kernels.\ \par 
But  $F$  does not depend on  $w$  and  $\displaystyle \ {\overline{\tilde
 S}}((z,w),\ \tilde a)$  is anti-holomorphic in  $w$  for  $z$
  fixed in  $\Omega ,$  so\ \par 
\quad \quad \quad \quad \quad   \[\displaystyle \ \frac{1}{\eta }\int_{\lbrace w\in {\mathbb{C}}^{k}::-\eta
 -r(z)\leq \left\vert{w}\right\vert ^{2}<-r(z)\rbrace }{{\overline{\tilde
 S}}((z,w),\ {\tilde a})\,dm(w)}\rightarrow {\overline{\tilde
 S((z,0),\ {\tilde a})}}v_{k}k(-r(z))^{k-1},\] \ \par 
by the proof of the subordination lemma, hence\ \par 
 \[\displaystyle \ \int_{\Omega }{f(z)\bar B_{k-1}(z,a)\,dm_{k-1}(z)}=\int_{\Omega
 }{f(z){\overline{\tilde S}}((z,0),\ {\tilde a})v_{k}k(-r(z))^{k-1}\,dm(z)}=\int_{\Omega
 }{f(z){\overline{\tilde S}}((z,0),\ {\tilde a})\,dm_{k-1}(z)}.\] \ \par 
\quad  So we have\ \par 
 \[\displaystyle \forall f\in A(\Omega ),\ \int_{\Omega }{f(z)({\overline{\tilde
 S}}((z,0),\ {\tilde a})-\bar B_{k-1}(z,a))\,dm_{k-1}(z)}=0,\] \ \par 
hence  $\displaystyle \tilde S((z,0),\ {\tilde a})-B_{k-1}(z,a)\perp
 A(\Omega )$  in  $A_{k-1}^{2}(\Omega ).$  But  $\displaystyle
 \tilde S((z,0),\ {\tilde a})-B_{k}(z,a)$  is holomorphic in  $z,$  hence\ \par 
\quad \quad \quad  $\displaystyle \forall z\in \Omega ,\ \tilde S((z,0),\ {\tilde
 a})=B_{k-1}(z,a).$ \ \par 
The second part is a direct application of the first part in
 the subordination lemma~\ref{subPrinCarl24}.  $\blacksquare $ \ \par 

\subsection{Interpolating sequences.}
\ \par 
\quad  For  $a\in \Omega ,$  let  $k_{a}(z):=S(z,a)$  denotes the Szeg\"o
 kernel of  $\Omega $  at the point  $a.$  It is also the reproducing
 kernel for  $H^{2}(\Omega ),$  i.e.\ \par 
\quad \quad \quad  $\displaystyle \forall a\in \Omega ,\ \forall f\in H^{2}(\Omega
 ),\ f(a)=\int_{\partial \Omega }{f(z)\bar k_{a}(z)\,d\sigma
 (z)}={\left\langle{f,\ k_{a}}\right\rangle}.$ \ \par 
\quad  Set  $\ {\left\Vert{k_{a}}\right\Vert}_{p}:={\left\Vert{k_{a}}\right\Vert}_{H^{2}(\Omega
 )}$  and:\ \par 
\begin{Definition}
 We say that the sequence  $\Lambda $  of points in  $\Omega
 $  is  $H^{p}(\Omega )$  interpolating if\par 
\quad \quad \quad  $\forall \lambda \in \ell ^{p}(\Lambda ),\ \exists f\in H^{p}(\Omega
 )::\forall a\in \Lambda ,\ f(a)=\lambda _{a}{\left\Vert{k_{a}}\right\Vert}_{p'},$
 \par 
with  $p'$  the conjugate exponent for  $\displaystyle p,\ \frac{1}{p}+\frac{1}{p'}=1.$
 \par 
\end{Definition}
\quad  We say that  $\Lambda $  has the linear extension property if
  $\Lambda $  is  $H^{p}(\Omega )$  interpolating and if moreover
 there is a bounded linear operator  $E\ :\ \ell ^{p}(\Lambda
 )\rightarrow H^{p}(\Omega )$   making the interpolation, i.e.\ \par 
\quad \quad \quad \quad \quad 	 $\forall \lambda \in \ell ^{p}(\Lambda ),\ E(\lambda )\in H^{p}(\Omega
 ),\ \forall a\in \Lambda ,\ E(\lambda )(a)=\lambda _{a}{\left\Vert{k_{a}}\right\Vert}_{p'}.$
 \ \par 
\quad  A weaker notion is the dual boundedness:\ \par 
\begin{Definition}
 We shall say that the sequence  $\Lambda $  of points in  $\Omega
 $  is dual bounded in  $H^{p}(\Omega )$  if there is a bounded
 sequence of elements in  $H^{p}(\Omega ),\ \lbrace \rho _{a}\rbrace
 _{a\in \Lambda }\subset H^{p}(\Omega )$  which dualizes the
 associated sequence of reproducing kernels, i.e.\par 
\quad \quad \quad  $\displaystyle \exists C>0::\forall a\in \Lambda ,\ {\left\Vert{\rho
 _{a}}\right\Vert}_{p}\leq C,\ \forall a,c\in \Lambda ,\ {\left\langle{\rho
 _{a},\ k_{c}}\right\rangle}=\delta _{a,c}{\left\Vert{k_{c}}\right\Vert}_{p'}.$
 \par 
\end{Definition}
\quad  If  $\Lambda $  is  $H^{p}(\Omega )$  interpolating then it
 is dual bounded in  $H^{p}(\Omega )\ :$  just interpolate the
 elements of the basic sequence in  $\ell ^{p}(\Lambda ).$ \ \par 
\quad  The converse is the crux of the characterization by Carleson~\cite{CarlInt58}
 of  $H^{\infty }({\mathbb{D}})$  interpolating sequences and
 the same by Shapiro \&  Shields~\cite{ShapShields61} for  $H^{p}({\mathbb{D}})$
  interpolating sequences  in  ${\mathbb{D}}.$ \ \par 
\quad \quad  	We do the same for the Bergman spaces.\ \par 
For  $k\in {\mathbb{N}}$  and  $a\in \Omega ,$  let  $b_{k,a}(z):=B_{k}(z,a)$
  denotes the Bergman kernel of  $\Omega $  at the point  $a.$
  It is also the reproducing kernel for  $A_{k}^{2}(\Omega ),$  i.e.\ \par 
\quad \quad \quad  $\displaystyle \forall a\in \Omega ,\ \forall f\in A_{k}^{2}(\Omega
 ),\ f(a)=\int_{\Omega }{f(z)\bar b_{k,\ a}(z)\,dm_{k}(z)}={\left\langle{f,\
 b_{k,\ a}}\right\rangle}.$ \ \par 
\quad  Now we set  $\ {\left\Vert{b_{k,\ a}}\right\Vert}_{p}:={\left\Vert{b_{k,\
 a}}\right\Vert}_{A_{k}^{p}(\Omega )}$  and:\ \par 
\begin{Definition}
 We say that the sequence  $\Lambda $  of points in  $\Omega
 $  is  $A_{k}^{p}(\Omega )$  interpolating if\par 
\quad \quad \quad  $\forall \lambda \in \ell ^{p}(\Lambda ),\ \exists f\in A_{k}^{p}(\Omega
 )::\forall a\in \Lambda ,\ f(a)=\lambda _{a}{\left\Vert{b_{k,a}}\right\Vert}_{p'},$
 \par 
with  $p'$  the conjugate exponent for  $\displaystyle p,\ \frac{1}{p}+\frac{1}{p'}=1.$
 \par 
\end{Definition}
\quad  We say that  $\Lambda $  has the linear extension property if
  $\Lambda $  is  $A_{k}^{p}(\Omega )$  interpolating and if
 moreover there is a bounded linear operator  $E\ :\ \ell ^{p}(\Lambda
 )\rightarrow A_{k}^{p}(\Omega )$   making the interpolation.\ \par 
\begin{Definition}
 We shall say that the sequence  $\Lambda $  of points in  $\Omega
 $  is dual bounded in  $A_{k}^{p}(\Omega )$  if there is a bounded
 sequence of elements in  $A_{k}^{p}(\Omega ),\ \lbrace \rho
 _{a}\rbrace _{a\in \Lambda }\subset A_{k}^{p}(\Omega )$  which
 dualizes the associated sequence of reproducing kernels, i.e.\par 
\quad \quad \quad  $\displaystyle \exists C>0::\forall a\in \Lambda ,\ {\left\Vert{\rho
 _{a}}\right\Vert}_{p}\leq C,\ \forall a,c\in \Lambda ,\ {\left\langle{\rho
 _{a},\ b_{k,c}}\right\rangle}=\delta _{a,c}{\left\Vert{b_{k,a}}\right\Vert}_{p'}.$
 \par 
\end{Definition}
\quad  Again if  $\Lambda $  is  $A_{k}^{p}(\Omega )$  interpolating
 then it is dual bounded in  $A_{k}^{p}(\Omega ):$  just interpolate
 the elements of the basic sequence in  $\ell ^{p}(\Lambda ).$ \ \par 

\subsection{Case of the unit disc  ${\mathbb{D}}$  in  ${\mathbb{C}}.$ }
\ \par 
\quad  In that case the interpolating sequences for  $H^{\infty }({\mathbb{D}})$
  where characterized by Carleson~\cite{CarlInt58} and for  $H^{p}({\mathbb{D}})$
  by Shapiro \&  Shields~\cite{ShapShields61}. The interpolating
 sequences for the Bergman spaces  $A_{k}^{p}({\mathbb{D}})$
  were characterized by Seip~\cite{Seip93}.\ \par 
In these cases it appears that dual boundedness implies interpolation.
 For Hardy spaces dual boundedness is easily seen to be equivalent
 to the Carleson condition and for Bergman spaces, it is proved
 by Schuster \&  Seip~\cite{SchuSeip98}.\ \par 

\subsection{General case.}
\ \par 
\quad  We shall apply the subordination lemma to interpolating sequences
 in general domains  $\Omega .$ \ \par 
Let  $\tilde \Omega $  be the lifted domain in  ${\mathbb{C}}^{n+k}$
  associated to  $\Omega .$  Let  $\tilde \Lambda $  be the sequence
  $\Lambda $  viewed in  $\tilde \Omega ,\ \tilde \Lambda :=\Lambda
 \subset \Omega \subset \tilde \Omega .$  Let us denote by  $k_{{\tilde
 a}}(z,w):=S((z,w),\ {\tilde a})$  the Szeg\"o kernel of  $\tilde
 \Omega ,\ $ for  ${\tilde a}=(a,0).$ \ \par 
\ \par 
\begin{Theorem}
~\label{subPrinFin542}Let  $\Omega $  be a domain in  ${\mathbb{C}}^{n}$
  and  $\tilde \Omega $  its lift to  ${\mathbb{C}}^{n+k}.$ 
 If  $\Lambda \subset \Omega $  is a sequence of points in  $\Omega
 ,$  let  $\tilde \Lambda $  be the sequence  $\Lambda $  viewed
 in  $\tilde \Omega ,\ \tilde \Lambda :=\Lambda \subset \Omega
 \subset \tilde \Omega .$  We have :\par 
\quad \quad  	(i)   $\Lambda $  is dual bounded in  $\displaystyle A_{k-1}^{p}(\Omega
 )$  iff   $\tilde \Lambda $  is dual bounded  in  $H^{p}(\tilde
 \Omega ).$ \par 
\quad \quad  	(ii)  $\Lambda $  is  $A_{k-1}^{p}(\Omega )$  interpolating
 iff   $\tilde \Lambda $  is  $H^{p}(\tilde \Omega )$  interpolating.\par 
\quad \quad  	(iii)  $\Lambda $  has the linear extension property in  $\displaystyle
 A_{k-1}^{p}(\Omega )$  iff  $\tilde \Lambda $   has the linear
 extension property in  $H^{p}(\tilde \Omega ).$ \par 
\end{Theorem}
\quad \quad  	Proof.\ \par 
\quad  For the {\sl (i) :} 	suppose that  $\Lambda $  is dual bounded
  $\displaystyle A_{k-1}^{p}(\Omega )$  and let  $\lbrace \rho
 _{a}\rbrace _{a\in \Lambda }\subset A_{k-1}^{p}(\Omega )\ $
  be the dual sequence to the sequence  $\lbrace b_{k-1,a}\rbrace
 _{a\in \Lambda };$  extend it to  $\tilde \Omega $  :\ \par 
\quad \quad \quad  $\forall a\in \Lambda ,\ \Gamma _{a}(z,w):=\rho _{a}(z),$ \ \par 
then the subordination lemma gives us that  $\ {\left\Vert{\Gamma
 _{a}}\right\Vert}_{H^{p}(\tilde \Omega )}\simeq {\left\Vert{\rho
 _{a}}\right\Vert}_{A_{k-1}^{p}(\Omega )}$  and we have, using
 corollary~\ref{subPrinCarl23},\ \par 
\quad \quad \quad  $\forall a,c\in \Lambda ,\ {\left\langle{\Gamma _{a},\ k_{\tilde
 c}}\right\rangle}={\left\langle{\Gamma _{a},\ S((\cdot ,\ 0),\
 \tilde c)}\right\rangle}={\left\langle{\rho _{a},\ B(\cdot ,\
 c)}\right\rangle}={\left\langle{\rho _{a},b_{k-1,c}}\right\rangle}=\delta
 _{ab}{\left\Vert{b_{k-1,c}}\right\Vert}_{p'},$ \ \par 
Because  $\Lambda $  is dual bounded in  $\displaystyle A_{k-1}^{p}(\Omega
 ).$  Then we have, by corollary~\ref{subPrinCarl23},\ \par 
\begin{equation}  \forall \tilde c=(c,0),\ c\in \Omega ,\ {\left\Vert{b_{k-1,c}}\right\Vert}_{A_{k-1}^{p}(\Omega
 )}\simeq {\left\Vert{k_{\tilde c}}\right\Vert}_{H^{p}(\Omega
 )},\label{subPrinFin441}\end{equation}\ \par 
hence\ \par 
\quad \quad \quad \quad \quad 	 $\displaystyle \forall a,c\in \Lambda ,\ {\left\langle{\Gamma
 _{a},\ k_{\tilde c}}\right\rangle}=\delta _{ac}{\left\Vert{b_{k-1,c}}\right\Vert}_{p'}\simeq
 \delta _{ac}{\left\Vert{k_{\tilde c}}\right\Vert}_{p'}$ \ \par 
hence  $\tilde \Lambda $  is dual bounded in  $H^{p}(\tilde \Omega ).$ \ \par 
\quad \quad  	Because we used only equivalences in this proof, it works also
 for the converse, hence if  $\tilde \Lambda $  is dual bounded
 in  $H^{p}(\tilde \Omega )$  then  $\Lambda $  is dual bounded
 in  $\displaystyle A_{k-1}^{p}(\Omega ).$ \ \par 
\quad \quad  	For the {\sl (ii) : }suppose that  $\tilde \Lambda $  is interpolating
 in  $H^{p}(\tilde \Omega ).$  We want to show that  $\Lambda
 $  is  $\displaystyle A_{k-1}^{p}(\Omega )$  interpolating,
 so let  $\mu =\lbrace \mu _{a}\rbrace _{a\in \Lambda }\in \ell
 ^{p}(\Lambda )$  the sequence to be interpolated. Set\ \par 
\quad \quad \quad \quad   $\lambda =\lbrace \lambda _{\tilde a}\rbrace _{a\in \Lambda
 }$  with  $\displaystyle \forall \tilde a\in \tilde \Lambda
 ,\ \lambda _{\tilde a}:=\mu _{a}{\times}\frac{{\left\Vert{b_{k-1,\
 a}}\right\Vert}_{A_{k-1}^{p'}(\Omega )}}{{\left\Vert{k_{\tilde
 a}}\right\Vert}_{H^{p'}(\tilde \Omega )}}\ ;$ \ \par 
then  $\lambda \in \ell ^{p}(\tilde \Lambda ),\ {\left\Vert{\lambda
 }\right\Vert}_{p}\simeq {\left\Vert{\mu }\right\Vert}_{p}$ 
 by~(\ref{subPrinFin441}).\ \par 
\quad  Let  $F\in H^{p}(\tilde \Omega )$  be the function making the
 interpolation of the sequence  $\lambda ,$  which exists because
  $\tilde \Lambda $  is  $H^{p}(\tilde \Omega )$  interpolating.
 It means that\ \par 
\quad \quad \quad \quad \quad 	 $F(\tilde a)=\lambda _{a}{\left\Vert{k_{\tilde a}}\right\Vert}_{H^{p'}(\tilde
 \Omega )}=\mu _{a}\ {\left\Vert{b_{k-1,\ a}}\right\Vert}_{A_{k-1}^{p'}(\Omega
 )}.$ \ \par 
Set  $\forall z\in \Omega ,\ f(z):=F(z,0)$  then we have\ \par 
\quad \quad \quad \quad \quad   $\forall a\in \Lambda ,\ f(a)=F(a,0)=F(\tilde a)=\mu _{a}\
 {\left\Vert{b_{k-1,\ a}}\right\Vert}_{A_{k-1}^{p'}(\Omega )}.$ \ \par 
Hence  $\Lambda $  is  $\displaystyle A_{k-1}^{p}(\Omega )$
  interpolating.\ \par 
\quad \quad  	Again the converse is straightforward because we use only equivalences.\ \par 
\quad  For the {\sl (iii)} : suppose that  $\tilde \Lambda $  has the
 bounded extension linear property, i.e. there is a linear operator
  $\tilde E\ :\ \ell ^{p}(\tilde \Lambda )\rightarrow H^{p}(\tilde
 \Omega )$  such that  $F(z,w):=\tilde E(\lambda )(z,w),$ \ \par 
\quad \quad \quad  $F\in H^{p}(\tilde \Omega ),\ \forall a\in \Lambda ,\ F({\tilde
 a})=\mu _{a}{\left\Vert{k_{{\tilde a}}}\right\Vert}_{H^{p'}(\tilde
 \Omega )},\ {\left\Vert{F}\right\Vert}_{H^{p}(\tilde \Omega
 )}\lesssim {\left\Vert{\mu }\right\Vert}_{p}.$ \ \par 
\quad  With the same notations  $\lambda $  and  $\mu $  as above,
 set  $f(z):=F(z,0)=\tilde E(\mu )(z,0)=:\ E(\lambda )(z),$ 
 then clearly  $E$  is linear in  $\lambda $  and then still
 using the subordination lemma we have\ \par 
\quad \quad \quad  $\ {\left\Vert{f}\right\Vert}_{A_{k-1}^{p}(\Omega )}\lesssim
 {\left\Vert{F}\right\Vert}_{H^{p}(\tilde \Omega )}$   $\lesssim
 {\left\Vert{\mu }\right\Vert}_{p}\simeq {\left\Vert{\lambda
 }\right\Vert}_{p}$ \ \par 
and\ \par 
\quad \quad \quad  $\forall a\in \Lambda ,\ f(a)=\mu _{\tilde a}{\left\Vert{k_{{\tilde
 a}}}\right\Vert}_{H^{p'}(\tilde \Omega )}=\lambda _{a}{\left\Vert{b_{k-1,\
 a}}\right\Vert}_{A_{k-1}^{p'}(\Omega )}.$ \ \par 
Hence  $\lambda \rightarrow E(\lambda )$  is bounded from  $\ell
 ^{p}(\Lambda )$  in  $\displaystyle A_{k-1}^{p}(\Omega )$  and
  $\Lambda $  is  $\displaystyle A_{k-1}^{p}(\Omega )$  interpolating
 with the linear extension.\ \par 
\quad \quad  	Again the converse is straightforward.  $\blacksquare $ \ \par 

\subsection{Application to strictly pseudo-convex domains.}
\ \par 
\quad \quad  	In~\cite{AmarExtInt06} we proved a general theorem on interpolating
 sequences in the spectrum of a uniform algebra. In the case
 of strictly pseudo-convex domains, it says that :\ \par 
\begin{Theorem}
 ~\label{subPrinFin543}If  $\Omega $  is a strictly pseudo-convex
 domain in  ${\mathbb{C}}^{n}$  and if  $\Lambda \subset \Omega
 $  is a dual bounded sequence of points in  $\displaystyle H^{p}(\Omega
 ),$  then, for any  $q<p,\ \Lambda $  is  $\displaystyle H^{q}(\Omega
 )$  interpolating with the linear extension property, provided
 that  $p=\infty $  or  $p\leq 2.$ \par 
\end{Theorem}
\quad \quad  	We have, as a consequence of the subordination lemma the following
 theorem.\ \par 
\begin{Theorem}
~\label{subPrinFin544}Let  $\Omega $  be a strictly pseudo-convex
 domain in  ${\mathbb{C}}^{n}$  and  $\Lambda \subset \Omega
 $  be a dual bounded sequence of points in  $\displaystyle A_{k}^{p}(\Omega
 ),$  then, for any  $q<p,\ \Lambda $  is  $\displaystyle A_{k}^{p}(\Omega
 )$  interpolating with the linear extension property, provided
 that  $p=\infty $  or  $p\leq 2.$ \par 
\end{Theorem}
\quad \quad  	Proof.\ \par 
Let  $\tilde \Omega $  be the lift of  $\Omega $  in  ${\mathbb{C}}^{n+k+1}$
  and  $\tilde \Lambda \subset \tilde \Omega $  the sequence
  $\Lambda $  viewed in  $\tilde \Omega .$  We apply theorem~\ref{subPrinFin542}
 {\sl (i)} to have that  $\tilde \Lambda $  is dual bounded in
  $H^{p}(\tilde \Omega )$  because  $\Lambda $  is dual bounded
 in  $\displaystyle A_{k}^{p}(\Omega ).$  Now we apply theorem~\ref{subPrinFin543}
 to get that  $\tilde \Lambda $  is  $H^{q}(\tilde \Omega )$
  interpolating with  $q<p,$  and has the bounded linear extension
 property, provided that  $p=\infty $  or  $p\leq 2.$  Then again
 theorem~\ref{subPrinFin542} {\sl (iii)} to get the same for
  $\Lambda $  in  $\displaystyle A_{k}^{p}(\Omega ).$   $\blacksquare $ \ \par 
\quad \quad  	We have a better result for the unit ball in  ${\mathbb{C}}^{n}\
 :$   in~\cite{intBall09} we proved\ \par 
\begin{Theorem}
~\label{subPrinFin545}If  $\Lambda $  is a dual bounded sequence
 in the unit ball  ${\mathbb{B}}$  of  ${\mathbb{C}}^{n}$  for
 the Hardy space  $H^{p}({\mathbb{B}}),$  then for any  $q<p,\
 \Lambda $  is  $H^{q}({\mathbb{B}})$  interpolating with the
 bounded linear extension property.\par 
\end{Theorem}
So copying the proof of theorem~\ref{subPrinFin544},  just replacing
 theorem~\ref{subPrinFin543} by theorem~\ref{subPrinFin545} we get :\ \par 
\begin{Theorem}
Let  $\Lambda $  be a dual bounded sequence in the unit ball
  ${\mathbb{B}}$  of  ${\mathbb{C}}^{n}$  for the Bergman space
  $A_{k}^{p}({\mathbb{B}}),$  then  for any  $q<p,\ S$  is  $A_{k}^{q}({\mathbb{B}})$
  interpolating with the bounded linear extension property.\par 
\end{Theorem}
\ \par 
\begin{Remark}
{\large   }If we apply this theorem in the unit disc  ${\mathbb{D}}$
  of   ${\mathbb{C}}$  we get that if  $\Lambda $  is  a dual
 bounded sequence in  $A_{k}^{p}({\mathbb{D}})$  then it is interpolating
 in  $A_{k}^{q}({\mathbb{D}})$  for any  $q<p.$  In this particular
 case, one variable, the Schuster-Seip theorem~\cite{SchuSeip98}
 says that we have the interpolation up to  $q=p.$ \par 
\end{Remark}
\ \par 

\subsection{Application to convex domains of finite type.}
\ \par 
\quad \quad  	To apply the general theorem on interpolating sequences in
 the spectrum of a uniform algebra to the case of convex domains
 of finite type in  ${\mathbb{C}}^{n},$  we need to have a precise
  knowledge of the good family of polydiscs associated to the
 domain and in~\cite{AmAspc09}, we proved\ \par 
\begin{Theorem}
 ~\label{subPrinFin546}If  $\Omega $  is a convex domain of finite
 type in  ${\mathbb{C}}^{n}$  and if  $\Lambda \subset \Omega
 $  is a dual bounded sequence of points in  $\displaystyle H^{p}(\Omega
 ),$  then, for any  $q<p,\ \Lambda $  is  $\displaystyle H^{q}(\Omega
 )$  interpolating with the linear extension property, provided
 that  $p=\infty $  or  $p\leq 2.$ \par 
\end{Theorem}
\quad  Then, again, copying the proof of theorem~\ref{subPrinFin544},
  just replacing theorem~\ref{subPrinFin543}  by theorem~\ref{subPrinFin546}
 we get :\ \par 
\begin{Theorem}
 If  $\Omega $  is a convex domain of finite type in  ${\mathbb{C}}^{n}$
  and if  $\Lambda \subset \Omega $  is a dual bounded sequence
 of points in  $A_{k-1}^{p}(\Omega )$  then, for any  $q<p,\
 \Lambda $  is  $A_{k-1}^{p}(\Omega )$  interpolating with the
 linear extension property, provided that  $p=\infty $  or  $p\leq 2.$ \par 
\end{Theorem}
\ \par 
\begin{Remark}
 We applied the subordination principle since 1978~\cite{DenAm78},
 ~\cite{amBerg78}  essentially in this case. For instance in~\cite{DenAm78}
 we used it to show that the interpolating sequences for  $H^{p}({\mathbb{B}}),$
  with  ${\mathbb{B}}$  the unit ball in  ${\mathbb{C}}^{n},\
 n\geq 2,$  are different for different values of  $p,$  opposite
 to the one variable case of  $H^{p}({\mathbb{D}}).$ \par 
\end{Remark}

\section{The  $H^{p}$ -Corona theorem for Bergman spaces.~\label{subPrinFin437}}

      Let  $\Omega $  be a domain in  ${\mathbb{C}}^{n}.$  We
 say that the  $H^{p}$ -Corona theorem is true for  $\Omega $
  if we have :\ \par 
\quad \quad \quad  $\forall g_{1},...,\ g_{k}\in H^{\infty }(\Omega )::\forall
 z\in \Omega ,\ \sum_{j=1}^{m}{\left\vert{g_{j}(z)}\right\vert
 }\geq \delta >0$ \ \par 
then\ \par 
\quad \quad \quad  $\forall f\in H^{p}(\Omega ),\ \exists (f_{1},...,\ f_{m})\in
 (H^{p}(\Omega ))^{m}::f=\sum_{j=1}^{m}{f_{j}g_{j}}.$ \ \par 
In the same vein, we say that the  $\displaystyle A_{k-1}^{p}(\Omega
 )$ -Corona theorem is true for  $\Omega $  if we have :\ \par 
\begin{equation}  \forall g_{1},...,\ g_{m}\in H^{\infty }(\Omega
 )::\forall z\in \Omega ,\ \sum_{j=1}^{m}{\left\vert{g_{j}(z)}\right\vert
 }\geq \delta >0\label{subPrinCarl25}\end{equation}\ \par 
then\ \par 
\quad \quad \quad  $\forall f\in A_{k-1}^{p}(\Omega ),\ \exists (f_{1},...,\ f_{m})\in
 (A_{k-1}^{p}(\Omega ))^{m}::f=\sum_{j=1}^{m}{f_{j}g_{j}}.$ \ \par 
We then have\ \par 
\begin{Theorem}
 Suppose that the  $H^{p}$ -Corona is true for the domain  $\tilde
 \Omega ,$  then the  $\displaystyle A_{k-1}^{p}(\Omega )$ -Corona
 theorem is also true for  $\Omega .$ \par 
\end{Theorem}
\quad  Proof.\ \par 
Let  $\tilde \Omega $  be the lifted domain ; then set\ \par 
\quad \quad \quad \quad  $\forall j=1,...,m,\ g_{j}\in H^{\infty }(\Omega ),\ f\in H^{p}(\Omega
 ),\ G_{j}(z,w):=g_{j}(z),\ F(z,w):=f(z).$ \ \par 
Clearly the  $G_{j}$  are in  $H^{\infty }(\tilde \Omega )$
  and by the subordination lemma,  $F\in H^{p}(\tilde \Omega
 ).$  Moreover, if the condition~(\ref{subPrinCarl25}) is true,
 we have  $\displaystyle \forall (z,w)\in \tilde \Omega ,\ \sum_{j=1}^{m}{\left\vert{G_{j}(z,w)}\right\vert
 }\geq \delta $  with the same  $\delta .$  So we can apply the
 hypothesis :\ \par 
\quad \quad \quad  $\exists (F_{1},...,\ F_{m})\in (H^{p}(\tilde \Omega ))^{m}::F=\sum_{j=1}^{m}{F_{j}G_{j}}.$
 \ \par 
Now set  $f_{j}(z)=F_{j}(z,0)$  then applying again the subordination
 lemma, we have\ \par 
\quad \quad \quad \quad \quad 	 $f(z)=F(z,0)=\sum_{j=1}^{m}{F_{j}(z,0)G_{j}(z,0)}=\sum_{j=1}^{m}{f_{j}(z)g_{j}(z).}$
  ???\ \par 
\ \par 

\subsection{Application to pseudo-convex domains.}
\ \par 
\begin{Corollary}
 We have the  $\displaystyle A_{k-1}^{p}(\Omega )$ -Corona theorem
  in the following cases :\par 
\quad  $\bullet $  with  $p=2$  if  $\Omega $  is a bounded weakly
 pseudo-convex domain in  ${\mathbb{C}}^{n};$ \par 
\quad  $\bullet $  with  $1<p<\infty $  if  $\Omega $  is a bounded
  strictly pseudo-convex domain in  ${\mathbb{C}}^{n}.$ \par 
\end{Corollary}
\quad  The first case because Andersson~\cite{Andersson94} (with a
 preprint in 1990) proved the  $H^{2}$  Corona theorem for  $\Omega
 $  bounded weakly pseudo-convex domain in  ${\mathbb{C}}^{n};$ \ \par 
\quad  the last one for two generators because we proved~\cite{amCor91}
 ( with~\cite{amWolff} already in 1980) the  $H^{p}$  Corona
 theorem for two generators in the ball ; for any number of generators
 because Andersson \&  Carlsson~\cite{AnderCarl00} (see also~\cite{AmMenDiv00})
 proved the  $H^{p}$  Corona theorem in this case.  $\blacksquare $ \ \par 

\section{Zeros set of the Nevanlinna-Bergman class.~\label{subPrinFin439}}
\quad  Let  $\Omega $  be a domain in  ${\mathbb{C}}^{n}$  and  $u$
  a holomorphic function in  $\Omega .$  Set  $X:=\lbrace z\in
 \Omega ::u(z)=0\rbrace $  the zero set of  $u$  and  $\Theta
 _{X}:=\partial \bar \partial \log \left\vert{u}\right\vert $
  its associated  $(1,1)$  current of integration.\ \par 
\begin{Definition}
 An analytic set   $X:=u^{-1}(0),\ u\in {\mathcal{H}}(\Omega
 ),$  in the domain  $\Omega $  is in the Blaschke class,  $X\in
 {\mathcal{B}}(\Omega ),$  if there is a constant  $C>0$  such that\par 
 \[\displaystyle \forall \beta \in \Lambda _{n-1,\ n-1}^{\infty
 }(\bar \Omega ),\ \left\vert{\int_{\Omega }{(-r(z))\Theta _{X}\wedge
 \beta }}\right\vert \leq C{\left\Vert{\beta }\right\Vert}_{\infty },\] \par 
where  $\displaystyle \Lambda _{n-1,\ n-1}^{\infty }(\bar \Omega
 )$  is the space of  $(n-1,n-1)$  continuous form in  $\bar
 \Omega ,$  equipped with the sup norm of the coefficients.\par 
\end{Definition}
\quad  If  $u\in {\mathcal{N}}(\Omega )$  then it is well known~\cite{zeroSkoda}
 that  $X$  is in the Blaschke class of  $\Omega .$ \ \par 
\quad  We do the analogue for the Bergman spaces :\ \par 
\begin{Definition}
 An analytic set  $\displaystyle X:=u^{-1}(0),\ u\in {\mathcal{H}}(\Omega
 ),$  in the domain  $\Omega $  is in the Bergman-Blaschke class,
  $X\in {\mathcal{B}}_{k-1}(\Omega ),$  if there is a constant
  $C>0$  such that\par 
 \[\displaystyle \forall \beta \in \Lambda _{n-1,\ n-1}^{\infty
 }(\bar \Omega ),\ \left\vert{\int_{\Omega }{(-r(z))^{k+1}\Theta
 _{X}\wedge \beta }}\right\vert \leq C{\left\Vert{\beta }\right\Vert}_{\infty
 },\] \par 
where  $\displaystyle \Lambda _{n-1,\ n-1}^{\infty }(\bar \Omega
 )$  is the space of  $(n-1,n-1)$  continuous form in  $\bar
 \Omega ,$  equipped with the sup norm of the coefficients.\par 
\end{Definition}
\quad  If  $u\in {\mathcal{N}}_{k-1}(\Omega )$  then  $X$  is in the
 Bergman-Blaschke class of  $\Omega ,$  for instance again by
 use the subordination lemma from the case  ${\mathcal{N}}(\tilde
 \Omega ).$ \ \par 
\quad  Hence exactly as for the Corona theorem we can set the definitions :\ \par 
we say that the {\sl Blaschke characterization is true for}
  $\Omega $  if we have :\ \par 
\quad \quad \quad  $X\in {\mathcal{B}}(\Omega )\Rightarrow \exists u\in {\mathcal{N}}(\Omega
 )$  such that  $X=\lbrace z\in \Omega ::u(z)=0\rbrace .$ \ \par 
And the same for the Bergman spaces :\ \par 
we say that the {\sl Bergman-Blaschke characterization is true
 for}  $\Omega $  if we have :\ \par 
\quad \quad \quad  $X\in {\mathcal{B}}_{k}(\Omega )\Rightarrow \exists u\in {\mathcal{N}}_{k}(\Omega
 )$  such that  $X=\lbrace z\in \Omega ::u(z)=0\rbrace .$ \ \par 
\begin{Theorem}
 Suppose that the Blaschke characterization is true for the lifted
 domain  $\tilde \Omega ,$  then the Bergman-Blaschke characterization
 is also true for  $\Omega .$ \par 
\end{Theorem}
\quad  Proof.\ \par 
Let  $\tilde \Omega $  be the lifted domain in  ${\mathbb{C}}^{n+k}$
  of  $\Omega $  ; then set   $X=u^{-1}(0),\ \Theta _{X}$  its
 associated current and suppose that  $X\in {\mathcal{B}}_{k}(\Omega ).$ \ \par 
\quad  This means that\ \par 
\quad \quad \quad  \[\displaystyle \forall \beta \in \Lambda _{n-1,\ n-1}^{\infty
 }(\bar \Omega ),\ \left\vert{\int_{\Omega }{(-r(z))^{k+1}\Theta
 _{X}\wedge \beta }}\right\vert \leq C{\left\Vert{\beta }\right\Vert}_{\infty
 }.\] \ \par 
Let\ \par 
\quad \quad \quad \quad \quad   $\forall w\in {\mathbb{C}}^{k},\ U(z,w):=u(z),\ \tilde X:=U^{-1}(0)\cap
 \tilde \Omega \subset \tilde \Omega ,\ \tilde \Theta _{\tilde
 X}=\partial \bar \partial \log  \left\vert{U}\right\vert \ ;$  \ \par 
we shall show that  $\tilde X\in {\mathcal{B}}(\tilde \Omega
 ).$  We have that  $\displaystyle \tilde \Theta _{\tilde X}$
  does not depend on  $w,$  hence\ \par 
 \[\displaystyle \forall \tilde \beta \in \Lambda _{n+k-1,\ n+k-1}^{\infty
 }({\overline{\tilde \Omega }}),\ A:=\int_{\tilde \Omega }{(-\tilde
 r(z,w))\tilde \Theta _{\tilde X}\wedge \tilde \beta }=\int_{\Omega
 }{\Theta _{X}(z)\wedge \int_{\left\vert{w}\right\vert ^{2}<-r(z)}{-(r(z)+\left\vert{w}\right\vert
 ^{2})\tilde \beta (z,w)}}.\] \ \par 
\quad  Because  $\Theta _{X}$  is a  $(1,1)$  current depending only
 on  $z,$  this means that in the integral in  $w$  we have only
 the terms containing  $dw_{1}\wedge d\bar w_{1}\wedge \cdot
 \cdot \cdot \wedge dw_{k}\wedge d\bar w_{k},$  the other terms
 being  $0$  against  $\Theta _{X}.$  	So this integral in  $w$
  gives a  $(n-1,n-1)$  form in  $z.$ \ \par 
\ \par 
\quad \quad  	Now set\ \par 
\quad \quad \quad \quad \quad 	 $\displaystyle \beta _{1}(z):=\int_{\left\vert{w}\right\vert
 ^{2}<-r(z)}{(1+\frac{\left\vert{w}\right\vert ^{2}}{-r(z)})\tilde
 \beta (z,w)},$ \ \par 
we have\ \par 
\quad \quad \quad \quad \quad 	 $\displaystyle A=\int_{\Omega }{\Theta _{X}(z)\wedge (-r(z))\beta
 _{1}(z)}$ \ \par 
and, because  $\displaystyle 1+\ \frac{\left\vert{w}\right\vert
 ^{2}}{-r(z)}<2$  in  $\lbrace \left\vert{w}\right\vert ^{2}<-r(z)\rbrace
 ,$  we have\ \par 
\quad \quad \quad \quad \quad 	 $\displaystyle \ \left\vert{\beta _{1}(z)}\right\vert \leq
 2{\left\Vert{\tilde \beta }\right\Vert}_{\infty }\int_{\left\vert{w}\right\vert
 ^{2}<-r(z)}{dm_{k}(w)\leq 2v_{k}{\left\Vert{\tilde \beta }\right\Vert}_{\infty
 }(-r(z))^{k},}$ \ \par 
because we get the volume in  ${\mathbb{C}}^{k}$  of the ball
 centered in  $0$  and of radius  $\ {\sqrt{-r(z)}}.$ \ \par 
\quad \quad  	Set  $\beta _{2}(z):=(-r(z))^{-k}\beta _{1}(z),$  we have 
 $\ {\left\Vert{\beta _{2}}\right\Vert}_{\infty }\leq 2v_{k}{\left\Vert{\tilde
 \beta }\right\Vert}_{\infty }$  and\ \par 
\quad \quad \quad \quad \quad 	 $\displaystyle A=\int_{\Omega }{\Theta _{X}(z)\wedge (-r(z))\beta
 _{1}(z)}=\int_{\Omega }{\Theta _{X}(z)\wedge (-r(z))^{k+1}\beta
 _{2}(z)}.$ \ \par 
\ \par 
\quad \quad  	We can apply the hypothesis  $X\in {\mathcal{B}}_{k-1}(\Omega
 )$  to the integral  $A$  :\ \par 
\quad \quad \quad  $\ \left\vert{A}\right\vert \leq {\left\Vert{\beta _{2}}\right\Vert}_{\infty
 }\lesssim 2{\left\Vert{\tilde \beta }\right\Vert}_{\infty },$ \ \par 
hence  $\tilde X\in {\mathcal{B}}(\tilde \Omega ).$ \ \par 
\quad  Now we apply the hypothesis of the theorem,\ \par 
\quad \quad \quad \quad  $\exists V\in {\mathcal{N}}(\tilde \Omega )::\tilde X=V^{-1}(0),$ \ \par 
and clearly  $X=V^{-1}(0)\cap \lbrace w=0\rbrace ,$  because
 if  $z\in X$  then  $\forall w::\left\vert{w}\right\vert ^{2}<-r(z),\
 (z,w)\in \tilde X.$  Hence we set\ \par 
\quad \quad \quad  $v(z):=V(z,0)\in {\mathcal{N}}_{k-1}(\Omega ),$ \ \par 
by the subordination lemma, and we are done.  $\blacksquare $ \ \par 

\subsection{Application to pseudo-convex domains.}
\ \par 
\begin{Corollary}
 The Bergman-Blaschke characterization is true  in the following cases :\par 
\quad  $\bullet $  if  $\Omega $  is a  strictly pseudo-convex domain
 in  ${\mathbb{C}}^{n}\ ;$ \par 
\quad  $\bullet $  if  $\Omega $  is a convex domain of finite type
 in  ${\mathbb{C}}^{n}.$ \par 
\end{Corollary}
\quad \quad  	Proof.\ \par 
The first case is true by the famous theorem proved by Henkin~\cite{zeroHenkin}
 and Skoda~\cite{zeroSkoda} which says that the Blaschke characterization
 is true for strictly bounded pseudo-convex domain in  ${\mathbb{C}}^{n}.$
 \ \par 
The second one because the Blaschke characterization is true
 for convex domain of finite strict type by a theorem of Bruna-Charpentier-Dupain~\cite{BruChaDup98}
 generalized to all convex domains of finite type by Cumenge~\cite{Cum01}
 and Diederich \&  Mazzilli~\cite{DiedMazz01}.  $\blacksquare $ \ \par 

\bibliographystyle{C:/texlive/2012/texmf-dist/bibtex/bst/base/plain}

\end{document}